\documentclass[twoside]{amsart}

\usepackage{latexsym, amsmath, amssymb, amsfonts, amscd, amsthm, mathrsfs}

\newcommand{\Ce}{{\mathbb C}}
\renewcommand{\Re}{{\mathbb R}}

\newcommand{\Ze}{{\mathbb Z}}
\newcommand{\Ne}{{\mathbb N}}

\newcommand{\LCTVS}{{\rm LCTVS}}
\newcommand{\LCA}{{\rm LCA}}
\newcommand{\LCSA}{{\rm LC$*$A}}
\newcommand{\Lie}{{\rm Lie}}

\newcommand{\pa}{\parallel}

\newcommand{\ie}{{\it i.e. }}

\theoremstyle{plain}
\newtheorem{theorem}{Theorem}[section]
\newtheorem{lemma}[theorem]{Lemma}
\newtheorem{proposition}[theorem]{Proposition}
\newtheorem{corollary}[theorem]{Corollary}

\theoremstyle{definition}
\newtheorem{remark}[theorem]{Remark}

\newtheorem{notation}[theorem]{Notation}
\newtheorem{definition}[theorem]{Definition}

\newtheorem*{acknowledgments}{Acknowledgments}

\begin{document}

\title{$\theta$-deformations as Compact Quantum Metric Spaces}
\author{Hanfeng Li}

\address{Department of Mathematics \\
University of Toronto \\
Toronto, ON M5S 3G3, CANADA} \email{hli@fields.toronto.edu}
\date{February 21, 2005}


\begin{abstract}
Let $M$ be a compact spin manifold with a smooth action of the
$n$-torus. Connes and Landi constructed $\theta$-deformations
$M_{\theta}$ of $M$, parameterized by $n\times n$ real
skew-symmetric matrices $\theta$. The $M_{\theta}$'s together with
the canonical Dirac operator $(D, \mathcal{H})$ on $M$ are an
isospectral deformation of $M$. The Dirac operator $D$ defines a
Lipschitz seminorm on $C(M_{\theta})$, which defines a metric on
the state space of $C(M_{\theta})$. We show that when $M$ is
connected, this metric induces the weak-$*$ topology. This means
that $M_{\theta}$ is a compact quantum metric space in the sense
of Rieffel.
\end{abstract}

\maketitle

\section{Introduction}
\label{intro:sec}

In noncommutative geometry there are many examples of noncommutative
spaces deformed from commutative spaces. However, for many of them
the Hochschild dimension,
which corresponds to the commutative notion of
dimension, is different from that of the original commutative space.
For instance, the $C^*$-algebras of the standard Podle\'s quantum
$2$-spheres and of the quantum $4$-spheres of \cite{BCT00}
are isomorphic to each other, and their
Hochschild dimension is zero \cite{MNW91}.

In \cite{CL01} Connes and Landi introduced a one-parameter
deformation $S^4_{\theta}$ of the $4$-sphere with the property
that the Hochschild dimension of $S^4_{\theta}$ equals that of
$S^4$. They also considered general $\theta$-deformations, which
was studied further by Connes and Dubois-Violette in \cite{CD01}
(see also \cite{Sitarz01}). In general, the $\theta$-deformation
$M_{\theta}$  of a manifold $M$ equipped with a smooth action of
the $n$-torus $T^n$ is determined by defining the algebra of
smooth functions  $C^{\infty}(M_{\theta})$ as the invariant
subalgebra (under the diagonal action of $T^n$) of the algebra
$C^{\infty}(M\times
T_{\theta}):=C^{\infty}(M)\hat{\otimes}C^{\infty}(T_{\theta})$ of
smooth functions on $M\times T_{\theta}$; here   $\theta$ is a
real skew-symmetric $n\times n$ matrix and $T_{\theta}$ is the
corresponding noncommutative $n$-torus.   This construction is a
special case of the strict deformation quantization  constructed
in \cite{Rieffel93}. When $M$ is a compact spin manifold,
Connes and Landi showed that the canonical Dirac operator $(D,
\mathcal{H})$ on $M$ and a deformed anti-unitary operator
$J_{\theta}$ together gives a spectral triple for
$C^{\infty}(M_{\theta})$, fitting it into Connes' noncommutative
Riemannian geometry framework \cite{Connes94, Connes96}. In
\cite{CD01} Connes and Dubois-Violette also showed how
$\theta$-deformations lead to compact quantum groups which are
deformations of various classical groups (see also \cite[Section
4]{Varilly01}).

In this paper we investigate the metric aspect of
$\theta$-deformation. The study of metric spaces in noncommutative
setting was initiated by Connes in \cite{Connes89} in the
framework of his spectral triple. The main ingredient of a
spectral triple is a Dirac operator $D$.  On the one hand, it
captures the differential structure by setting $df=[D, f]$. On the
other hand, it enables us to recover the Lipschitz seminorm $L$,
which is usually defined as
\begin{eqnarray} \label{dist to Lip:eq}
L(f):=\sup \{\frac{|f(x)-f(y)|}{\rho(x, y)}:x \neq y\},
\end{eqnarray}
where $\rho$ is the geodesic metric on the Riemannian manifold,
instead by means of $L(f)=\pa [D, f]\pa$, and then one recovers the metric $\rho$  by
\begin{eqnarray} \label{Lip to dist:eq}
\rho(x, y)=\sup_{L(f)\le 1} |f(x)-f(y)|.
\end{eqnarray}
In \cite[Section 1]{Connes89} Connes went further by considering
the (possibly $+\infty$-valued) metric on the state space of the
algebra defined by (\ref{Lip to dist:eq}). Motivated by what
happens to ordinary compact metric spaces, in \cite{Rieffel98b,
Rieffel99b, Rieffel00} Rieffel introduced ``compact quantum metric
spaces'' (see Definition~\ref{C*QCM:def} below) which requires the
metric on the state space to induce the $w^*$-topology.
Many examples of compact quantum metric spaces have been
constructed, mostly from ergodic actions of compact groups
\cite{Rieffel98b} or group algebras \cite{Rieffel02, Rieffel03O}.
Usually it is quite difficult to find out whether a specific
seminorm $L$ on a unital $C^*$-algebra gives a quantum metric,
i.e., whether the metric defined by (\ref{Lip to dist:eq}) on the
state space induces the $w^*$-topology.

Denote by $L_{\theta}$ the seminorm on $C(M_{\theta})$ determined by the Dirac operator $D$
(see Definition~\ref{L_{theta}:def} below for detail).
Notice that when $M$ is connected the geodesic distance makes $M$ into
a metric space. Then our main theorem in this paper is:

\begin{theorem} \label{theta-deform metric:thm}
Let $M$ be a connected compact spin manifold with a smooth action
of $\mathbb{T}^n$. For every real skew-symmetric $n\times n$
matrix $\theta$ the pair $(C(M_{\theta}),L_{\theta})$ is a
$C^*$-algebraic compact quantum metric space.
\end{theorem}

Motivated by questions in string theory, Rieffel also introduced a
notion of quantum Gromov-Hausdorff distance for compact quantum
metric spaces \cite{Rieffel00, Rieffel01}. It has many nice
properties. Using the quantum Gromov-Hausdorff distance one can
discuss the continuity of $\theta$-deformations (with respect to
the parameter $\theta$) in a concrete way. This will be done in
\cite{Li10}.

This paper is organized as follows. We shall use heavily the
theory of locally convex topological vector spaces ({\LCTVS}). In
Section~\ref{Prelim:sec} we review some facts about {\LCTVS},
Clifford algebras, and Rieffel's theory of compact quantum metric
spaces. Connes and Dubois-Violette's formulation of
$\theta$-deformations is reviewed in
Section~\ref{ConsConnesDubois:sec}. In Section~\ref{Lip
action:sec} we prove a general theorem showing that in the
presence of a compact group action, sometimes we can reduce the
study of a given seminorm to its behavior on the isotypic
components of this group action.
Section~\ref{DiffSN:sec} contains the main part of our proof of
Theorem~\ref{theta-deform metric:thm}, where
we study various differential operators
to derive certain formulas. Finally, Theorem~\ref{theta-deform metric:thm} is proved in
Section~\ref{proof:sec}.

Throughout this paper $G$ will be a nontrivial compact group with
identity $e_G$, endowed with the normalized Haar measure.
Denote by $\hat{G}$ the dual of $G$, and by $\gamma_0$
the trivial representation.
For any $\gamma\in \hat{G}$ let $\chi_{\gamma}$ be the corresponding character on
$G$, and let $\bar{\gamma}$ be the contragradient representation .
 For any $\gamma
\in \hat{G}$ and any representation of $G$ on some complex vector space
$V$, we denote by $V_{\gamma}$ the $\gamma$-isotypic component of $V$. If
$\mathcal{J}$ is a finite subset of $\hat{G}$, we also let
$V_{\mathcal{J}}=\sum_{\gamma\in \mathcal{J}}V_{\gamma}$, and let
$\bar{\mathcal{J}}=\{\bar{\gamma}:\gamma \in \mathcal{J}\}$.

\begin{acknowledgments}
  This is part of my Ph.D. dissertation submitted to UC Berkeley in 2002.
  I am indebted to my advisor, Professor Marc Rieffel, for many helpful discussions, suggestions,
  and for his support throughout my time at Berkeley.
I also thank Thomas Hadfield and Fr\'ed\'eric Latr\'emoli\`ere
  for valuable conversations.
\end{acknowledgments}

\section{Preliminaries}
\label{Prelim:sec}

In this section we review some facts about locally convex
topological vector spaces (\LCTVS), Clifford algebras, and
Rieffel's theory of compact quantum metric spaces.

\subsection{Locally convex
topological vector spaces}\label{LCTVS:sub}

 We recall first some facts about
{\LCTVS}. The reader is referred to \cite[Chapters 5 and
43]{Treves67} for detailed information about completion and tensor
products of {\LCTVS}. Throughout this paper, our {\LCTVS}   will
all be Hausdorff.

For any {\LCTVS} $V$ and $W$, one can define the \emph{projective
tensor product} of $V$ and $W$, denoted by $V\otimes_{\pi}W$, as
the vector space $V\otimes W$ equipped with the so called projective
topology. $V\otimes_{\pi}W$ is also a {\LCTVS}, and one can form
the completion $V\hat{\otimes}_{\pi}W$.

For continuous linear
maps $\psi_j:V_j\rightarrow W_j$ ($j=1,2$) between
{\LCTVS}, the tensor product linear
map
$\psi_1\otimes_{\pi}
\psi_2:V_1\otimes_{\pi} V_2\rightarrow W_1\otimes_{\pi} W_2$ is also continuous and
extends to a continuous linear map $\psi_1\hat{\otimes}_{\pi}\psi_2:
V_1\hat{\otimes}_{\pi}V_2\rightarrow W_1\hat{\otimes}_{\pi}W_2$.

Let $V$ be a {\LCTVS}, and let $\alpha$ be an action
of a topological $G$ on $V$ by automorphisms. We say that the action $\alpha$
is \emph{continuous}
if the map $G\times V\rightarrow V$ given by $(x, v)\mapsto
\alpha_x(v)$ is (jointly) continuous. Let $V$ (resp. $W$) be a {\LCTVS} and $\alpha$
(resp. $\beta$) be a continuous action
of $G$ on $V$ (resp. $W$). Then the tensor product action $\alpha \hat{\otimes}_{\pi} \beta$
of $G$ on $V\hat{\otimes}_{\pi}W$ is easily seen to be
continuous.

A \emph{locally convex algebra} (\LCA) \cite{Connes86} is a
{\LCTVS} V with an algebra structure such that the multiplication
$V\times V\rightarrow V $ is (jointly) continuous. If furthermore
$V$ is a $*$-algebra and the $*$-operation $*:V\rightarrow V$ is
continuous, let us say that $V$ is a \emph{locally convex
$*$-algebra} (\LCSA). A \emph{locally convex left $V$-module} of
$V$ is a left $V$-module $W$ such that the action $V\times
W\rightarrow W$ is (jointly) continuous. For a smooth manifold
$M$, the space of (possibly unbounded) smooth functions
$C^{\infty}(M)$ equipped with usual Fr\'echet space topology is a
{\LCSA}. For a smooth vector bundle $E$ over $M$, the space of
smooth sections $C^{\infty}(M, E)$ is a locally convex
$C^{\infty}(M)$-bimodule. If furthermore $E$ is an algebra bundle
with fibre algebras being finite-dimensional, then $C^{\infty}(M,
E)$ is also a {\LCA}. We shall need Proposition~\ref{tensor prod
of alg:lemma} below.

\begin{lemma} \label{completion tensor prod:lemma}
Let $V$ and $W$ be two {\LCTVS}. Denote by $\hat{V}$ and $\hat{W}$
the completion of $V$ and $W$ respectively. Then
\begin{eqnarray*}
\hat{V}\hat{\otimes}_{\pi}\hat{W}=V\hat{\otimes}_{\pi}W.
\end{eqnarray*}
\end{lemma}
\begin{proof}
The natural linear maps $\iota_V:V\hookrightarrow \hat{V}$ and
$\iota_W:W\hookrightarrow \hat{W}$ are continuous, so we have the
continuous linear map
$\iota_V\hat{\otimes}_{\pi}\iota_W:V\hat{\otimes}_{\pi}W\rightarrow
\hat{V}\hat{\otimes}_{\pi}\hat{W}$, which is the unique continuous
extension of $\iota_V\otimes\iota_W:V\otimes W\rightarrow
\hat{V}\otimes \hat{W}$.

Let $v_0\in \hat{V}$ (resp. $w_0\in \hat{W}$) and a net
$\{v_j\}_{j\in I}$ (resp. $\{w_j\}_{j\in I}$) in $V$ (resp. $W$)
converging to $v_0$ (resp. $w_0$). Let $\mathfrak{p}$ (resp.
$\mathfrak{q}$) be a continuous seminorm on $V$ (resp. $W$).
Consider the continuous tensor product seminorm $
\mathfrak{p}\hat{\otimes}_{\pi}\mathfrak{q}$ on $V\hat{\otimes}_{\pi}W$
defined by
\begin{eqnarray*}
(\mathfrak{p}\hat{\otimes}_{\pi} \mathfrak{q})(\eta)=
\inf \sum_j\mathfrak{p}(v'_j)\mathfrak{q}(w'_j)
\end{eqnarray*}
for all $\eta\in V\otimes_{\pi} W$, where the infimum is taken over all
finite sets of pairs $(v'_k, w'_k)$ such that
\begin{eqnarray*}
\eta=\sum_kv'_k\otimes w'_k.
\end{eqnarray*}
It satisfies
\begin{eqnarray*}
(\mathfrak{p}\hat{\otimes}_{\pi} \mathfrak{q})(v\otimes
w)=\mathfrak{p}(v)\mathfrak{q}(w)
\end{eqnarray*}
for all $v\in V$ and $w\in W$ \cite[Proposition 43.1]{Treves67}.
In particular, we have
\begin{eqnarray*}
(\mathfrak{p}\hat{\otimes}_{\pi} \mathfrak{q})(v_j\otimes
w_j-v_{j'}\otimes w_{j'}) &=& (\mathfrak{p}\hat{\otimes}_{\pi}
\mathfrak{q})((v_j-v_{j'})\otimes
w_j+v_{j'}\otimes (w_j-w_{j'})) \\
&\le &
\mathfrak{p}(v_j-v_{j'})\mathfrak{q}(w_j)+\mathfrak{p}(v_{j'})\mathfrak{q}(w_j-w_{j'})\to  0
\end{eqnarray*}
as $j,j'\to \infty$.
Since such $\mathfrak{p}\hat{\otimes}_{\pi}\mathfrak{q}$ form a
basis of continuous seminorms on $V\hat{\otimes}_{\pi}W$
\cite[page 438]{Treves67}, the net $\{v_j\otimes w_j\}_{j\in I}$
is a Cauchy net in $V\hat{\otimes}_{\pi}W$. Then it converges to
some element in $V\hat{\otimes}_{\pi}W$. Let $\varphi(v_0,
w_0)=\lim_{j\to \infty} (v_j\otimes w_j)$. Clearly $\varphi(v_0,
w_0)$ doesn't depend on the choice of the nets $\{v_j\}_{j\in I}$
and $\{w_j\}_{j\in I}$. So the map $\varphi:\hat{V}\times
\hat{W}\rightarrow V\hat{\otimes}_{\pi}W$ is well-defined. It is
easy to see that $\varphi$ is bilinear and is an extension of the
natural map $V\times W\rightarrow V\hat{\otimes}_{\pi}W$. Denote
the extension of $\mathfrak{p}$ (resp. $\mathfrak{q}$) on
$\hat{V}$ (resp. $\hat{W}$) still by $\mathfrak{p}$ (resp.
$\mathfrak{q}$). Notice that
\begin{eqnarray*}
(\mathfrak{p}\hat{\otimes}_{\pi}\mathfrak{q})(\varphi(v_0, w_0))
&=&
(\mathfrak{p}\hat{\otimes}_{\pi}\mathfrak{q})(\lim_{j\to \infty}
(v_j\otimes w_j))
=\lim_{j\to
  \infty}(\mathfrak{p}\hat{\otimes}_{\pi}\mathfrak{q})(v_j\otimes
w_j)\\
&=& \lim_{j\to \infty}\mathfrak{p}(v_j)\mathfrak{q}(w_j)= \mathfrak{p}(v)\mathfrak{q}(w).
\end{eqnarray*}
So $\varphi$ is continuous, and hence the associated linear map
$\hat{V}\otimes_{\pi}\hat{W}\rightarrow V\hat{\otimes}_{\pi}W$ is
continuous \cite[Proposition 43.4]{Treves67}. Consequently, we
have the continuous extension
$\psi:\hat{V}\hat{\otimes}_{\pi}\hat{W}\rightarrow
V\hat{\otimes}_{\pi}W$ \cite[Theorem 5.2]{Treves67}.

Notice that $V\otimes W$ is dense in both
$\hat{V}\hat{\otimes}_{\pi}\hat{W}$
and $V\hat{\otimes}_{\pi}W$. Clearly $\psi$ and
$\iota_V\hat{\otimes}_{\pi}\iota_W$
are inverse to each other
when restricted to $V\otimes W$. It follows immediately
that $\psi$ and $\iota_V\hat{\otimes}_{\pi}\iota_W$ are isomorphisms
inverse to each other between
$\hat{V}\hat{\otimes}_{\pi}\hat{W}$
and $V\hat{\otimes}_{\pi}W$.
\end{proof}

\begin{lemma} \label{tensor prod of bilinear:lemma}
Let $V_j, W_j, H_j$ ($j=1,2$) be {\LCTVS}, and let
$\psi_j:V_j\times W_j\rightarrow H_j$ be continuous bilinear maps;
then the bilinear map
\begin{eqnarray*}
\psi_1\otimes \psi_2:(V_1\otimes
V_2)\times (W_1\otimes W_2)\rightarrow H_1\otimes H_2
\end{eqnarray*}
extends to a continuous bilinear map
\begin{eqnarray*}
\psi_1\hat{\otimes}_{\pi}\psi_2:(V_1\hat{\otimes}_{\pi}
V_2)\times (W_1\hat{\otimes}_{\pi} W_2)\rightarrow
H_1\hat{\otimes}_{\pi} H_2.
\end{eqnarray*}
\end{lemma}
\begin{proof}
We have the associated continuous linear map
$\varphi_j:V_j\otimes_{\pi}W_j\rightarrow H_j, j=1,2$
\cite[Proposition 43.4]{Treves67} and hence the continuous linear
map
\begin{eqnarray*}
\varphi_1\hat{\otimes}_{\pi}\varphi_2:(V_1\otimes_{\pi}W_1)
\hat{\otimes}_{\pi}(V_2\otimes_{\pi}W_2)\rightarrow
H_1\hat{\otimes}_{\pi} H_2.
\end{eqnarray*}
By the associativity of the projective tensor product and
 Lemma~\ref{completion tensor prod:lemma} we have
\begin{eqnarray*}
& &(V_1\otimes_{\pi}W_1)\hat{\otimes}_{\pi}(V_2\otimes_{\pi}W_2)\\
&=&
((V_1\otimes_{\pi}W_1)\otimes_{\pi}V_2)\hat{\otimes}_{\pi}W_2
= ((V_1\otimes_{\pi}V_2)\otimes_{\pi}W_1)\hat{\otimes}_{\pi}W_2\\
&=& (V_1\otimes_{\pi}V_2)\hat{\otimes}_{\pi}(W_1\otimes_{\pi}W_2)
= (V_1\hat{\otimes}_{\pi}V_2)\hat{\otimes}_{\pi}(W_1\hat{\otimes}_{\pi}W_2).
\end{eqnarray*}
So we get a continuous linear map
$(V_1\hat{\otimes}_{\pi}V_2)\hat{\otimes}_{\pi}
(W_1\hat{\otimes}_{\pi}W_2)\rightarrow H_1\hat{\otimes}_{\pi}H_2$,
which is equivalent to a continuous bilinear map
$(V_1\hat{\otimes}_{\pi}V_2)\times
(W_1\hat{\otimes}_{\pi}W_2)\rightarrow H_1\hat{\otimes}_{\pi}H_2$.
Clearly this extends the bilinear map $\psi_1\otimes
\psi_2:(V_1\otimes V_2)\times (W_1\otimes W_2)\rightarrow
H_1\otimes H_2$.
\end{proof}

\begin{proposition} \label{tensor prod of alg:lemma}
Let $V$ and $W$ be
{\LCA}. Then
$V\hat{\otimes}_{\pi}W$ is also a
{\LCA} extending the
natural algebra structure on $V\otimes W$. If both $V$ and $W$ are
{\LCSA}, so is $V\hat{\otimes}_{\pi}W$. If $H$ is a
locally convex left $V$-module, then
$H\hat{\otimes}_{\pi}W$ is a locally convex left
$V\hat{\otimes}_{\pi}W$-module.
\end{proposition}
\begin{proof}By Lemma~\ref{tensor prod of bilinear:lemma} we have the
  continuous bilinear map
\begin{eqnarray*}
(V\hat{\otimes}_{\pi}W)\times (V\hat{\otimes}_{\pi}W)\rightarrow
V\hat{\otimes}_{\pi}W
\end{eqnarray*}
extending the multiplication of $V\otimes W$. Since $V\otimes W$ is
dense in $V\hat{\otimes}_{\pi}W$, clearly the above bilinear map is
associative. In other words, $V\hat{\otimes}_{\pi}W$ is
 a {\LCA}. The assertion about modules can be proved
 in the same way.

If both $V$ and $W$ are {\LCSA},
then we have the tensor product of the
$*$-operations $V\hat{\otimes}_{\pi}W\rightarrow
V\hat{\otimes}_{\pi}W$. Since it extends the natural $*$-operation on
$V\otimes W$, it is easy to check that it is compatible with the algebra
structure. So $V\hat{\otimes}_{\pi}W$ is a {\LCSA}.
\end{proof}

For any {\LCTVS} $V$ and $W$, one can also define the
\emph{injective tensor product} $V\otimes_{\epsilon}W$ of $V$ and
$W$, and form the completion $V\hat{\otimes}_{\epsilon}W$. Let us
say that a continuous linear map $\psi:V\rightarrow W$ is an
isomorphism of $V$ into $W$ if $\psi$ is injective and
$\psi:V\rightarrow \psi(V)$ is a homeomorphism of topological
spaces. The only property about injective tensor product we shall
need is that if $\psi_j$ is an isomorphism of $V_j$ into $W_j$ for
$j=1,2$, then the corresponding tensor product linear map
$\psi_1\hat{\otimes}_{\epsilon}\psi_2$ is an isomorphism of
$V_1\hat{\otimes}_{\epsilon}V_2$  into
$W_1\hat{\otimes}_{\epsilon}W_2$ \cite[Proposition
43.7]{Treves67}.

Let $n\ge 2$, and
let $\theta$ be a real skew-symmetric $n\times n$ matrix. Denote
by $\mathcal{A}_{\theta}$ the corresponding  quantum torus
\cite{Rieffel88, Rieffel90}.  It could be described as follows.
Let $\omega_{\theta}$ denote the skew-symmetric bicharacter on $\Ze^n$
defined by
\begin{eqnarray*}
\omega_{\theta}(p, q)=e^{i\pi p\cdot \theta q}.
\end{eqnarray*}
For each $p\in \Ze^n$ there is a unitary $u_p$ in
$\mathcal{A}_{\theta}$. And $\mathcal{A}_{\theta}$ is generated by
these unitaries with the relation
\begin{eqnarray*}
u_pu_q=\omega_{\theta}(p, q)u_{p+q}.
\end{eqnarray*}
So one may think of vectors in $\mathcal{A}_{\theta}$ as some kind of
functions on $\Ze^n$ .  The
$n$-torus $\mathbb{T}^n$ has a canonical ergodic action
$\tau$ on $\mathcal{A}_{\theta}$. Notice that $\Ze^n$ is
the dual group of $\mathbb{T}^n$. We denote the duality by $\left<p,
x\right>$ for $x\in \mathbb{T}^n$ and $p\in \Ze^n$. Then
$\tau$ is determined by
\begin{eqnarray*}
\tau_x(u_p)=\left<p, x\right>u_p.
\end{eqnarray*}
The set $\mathcal{A}^{\infty}_{\theta}$ of smooth vectors for the
action $\tau$ is exactly the Schwarz space $\mathcal{S}(\Ze^n)$
\cite{Connes80}. Let $X_1, \cdots, X_n$ be a basis for the Lie
algebra of $\mathbb{T}^n$. Then we have the differential
$\partial_{X_j}(f)$ for each $f\in \mathcal{A}^{\infty}_{\theta}$
and $1\le j\le n$. For each $k\in \Ne$ define a seminorm,
$\mathfrak{q}_k$, on $\mathcal{A}^{\infty}_{\theta}$ by
\begin{eqnarray*}
\mathfrak{q}_k:=\max_{|\overrightarrow{m}|\le k}\pa
\partial^{m_1}_{X_1}\cdots \partial^{m_n}_{X_n}(f)\pa.
\end{eqnarray*}
Clearly $\mathcal{A}^{\infty}_{\theta}$ is a complete {\LCSA}
equipped with the topology defined by these $\mathfrak{q}_k$'s. On
the other hand, it is easy to see that this topology is the same
as the usual topology on $\mathcal{S}(\Ze^n)$. Thus
$\mathcal{A}^{\infty}_{\theta}$ is a \emph{nuclear} space
\cite[Theorem 51.5]{Treves67}, which means that for every {\LCTVS}
$V$ the injective and projective topologies on $V\otimes
\mathcal{A}^{\infty}_{\theta}$ coincide \cite[Theorem
50.1]{Treves67}. So we shall simply use $V\otimes
\mathcal{A}^{\infty}_{\theta}$ to denote the (projective or
injective) topological tensor product. The algebraic tensor
product will be denoted by $V\otimes_{alg}
\mathcal{A}^{\infty}_{\theta}$.

We shall need to integrate continuous functions with values in a
{\LCTVS}. For our purpose, it suffices to use the
Riemann integral. Though this should be well-known, we have not been able to
find any reference in the literature. So we include a definition here.

\begin{lemma}  \label{integ:lemma}
Let $X$ be a compact space with a probability measure $\mu$.
Let
\begin{eqnarray*}
I:=&\{&\{X_1, {\cdots}, X_k\}:X_1, {\cdots}, X_k \mbox{ are
disjoint
   measurable subsets of } X, \\
& &k\in \Ne,\,  \cup^k_{j=1}X_j=X\}
\end{eqnarray*}
be the set of all finite partitions of $X$ into measurable subsets
with the fine order, \ie
\begin{eqnarray*}
\{X_1, {\cdots}, X_k\}\ge \{X'_1, {\cdots}, X'_{k'}\} \mbox{ if
and only if every } X_j \mbox{ is contained in some } X'_{j'}.
\end{eqnarray*}
Let $V$ be a complete {\LCTVS}, and let $f:X\rightarrow V$ be a
continuous map. For each $\{X_1, {\cdots}, X_k\}$ in $I$ pick an
$x_j\in X_j$ for each $j$, and let
\begin{eqnarray*}
v_{\{X_1, {\cdots}, X_k\}}=\sum^k_{j=1}\mu(X_j)f(x_j).
\end{eqnarray*}
Then $\{v_{\{X_1, {\cdots}, X_k\}}\}_{\{X_1, {\cdots}, X_k\}\in
I}$ is a
  Cauchy net in $V$, and its limit doesn't depend on the choice of the
  representatives
$x_1, {\cdots}, x_k$.
\end{lemma}
\begin{proof} Let a continuous seminorm $\mathfrak{p}$ on $V$ and
an $\epsilon>0$ be given. For each $x\in X$ there is an open
neighborhood $\mathcal{U}_x$ of $x$ such that
$\mathfrak{p}(f(x)-f(y))\le \epsilon$ for all $y\in\mathcal{U}_x$.
Since $X$ is compact, we can cover $X$ with finitely many such
$\mathcal{U}_x$, say $\mathcal{U}_{x_1}, {\cdots},
\mathcal{U}_{x_k}$. Let $X_1=\mathcal{U}_{x_1}$ and
$X_j=\mathcal{U}_{x_j}\setminus \cup^{j-1}_{s=1}X_s$ inductively
for all $2\le j\le k$. Then $\{X_1, {\cdots}, X_k\}$ is a finite
partition of $X$. For any $\{X'_1, {\cdots}, X'_{k'}\}\ge \{X_1,
{\cdots}, X_k\}$, clearly $\mathfrak{p}(v_{\{X'_1, {\cdots},
X'_{k'}\}}-v_{\{X_1, {\cdots}, X_k\}})\le
  2\epsilon$ no matter how we choose the representatives for
$\{X'_1, {\cdots}, X'_{k'}\}$ and $\{X_1, {\cdots}, X_k\}$. This
gives the desired result.
\end{proof}

\begin{definition} \label{integ:def}
Let $X$ be  a compact space with a probability measure $\mu$, and let
$f$ be a
continuous function from $X$ into a complete
{\LCTVS} $V$. The integration of $f$ over $X$, denoted by
$\int_Xf\,d\mu$, is defined as the limit in Lemma~\ref{integ:lemma}.
\end{definition}

The next proposition is obvious:

\begin{proposition}  \label{linear of integ:prop}
Let $X$ be  a compact space with a probability measure $\mu$, and let $f_1,
f_2$ be
continuous functions from $X$ into a complete
{\LCTVS} $V$. Then
\begin{eqnarray*}
\int_X(f_1+f_2)\, d\mu&=&\int_Xf_1\,d\mu+\int_Xf_2\,d\mu, \\
\int_X\lambda f_1\,d\mu&=&\lambda\int_Xf_1\,d\mu
\end{eqnarray*}
for any scalar $\lambda$.
If $\psi:V\rightarrow W$ is a continuous linear map from $V$ into
another complete
{\LCTVS} $W$, then
\begin{eqnarray*}
\int_X\psi\circ f_1\, d\mu=\psi(\int_Xf_1\,d\mu).
\end{eqnarray*}
\end{proposition}

It is also easy to verify the analogue of the
fundamental theorem of calculus:

\begin{proposition} \label{integ to diff:prop}
Let $f$ be a continuous map from $[0,1]$ to a complete
{\LCTVS} $V$. Then
\begin{eqnarray*}
f(0)=\lim_{t\to 0}\frac{\int^t_0f(s)\, ds}{t}.
\end{eqnarray*}
\end{proposition}

\subsection{Clifford algebras}\label{clifford:sub}
Next we recall some facts about Clifford algebras \cite[Chapter
1]{GM91} \cite[Section 1.8]{Jost98}.

Let $V$ be a real vector space of dimension $m$ equipped with a
positive-definite inner product. The corresponding \emph{Clifford algebra},
denoted by $Cl(V)$, is the quotient of the tensor algebra
$\oplus_{k\ge 0}V\otimes \cdots \otimes V$ generated by $V$ by the two
sided ideal generated by all elements of the form $v\otimes v+\pa v\pa
^2$ for $v\in V$. The \emph{complexified Clifford algebra}, denoted by
$Cl^{\Ce}(V)$, is defined
as $Cl^{\Ce}(V):=Cl(V)\otimes_{\Re}\Ce$.

$Cl^{\Ce}(V)$ has a natural finite-dimensional $C^*$-algebra
structure \cite[Theorem 1.7.35]{GM91}. Denote by $SO(V)$ the group
of isometries of $V$ preserving the orientation. For each $g\in
SO(V)$ the isometry $g:V\rightarrow V$ induces an algebra
isomorphism $Cl(V)\rightarrow Cl(V)$ and a $C^*$-algebra
isomorphism $Cl^{\Ce}(V)\rightarrow Cl^{\Ce}(V)$. In this way
$SO(V)$ acts on $Cl(V)$ and $Cl^{\Ce}(V)$.

Recall that a state $\varphi$ on a $C^*$-algebra $\mathcal{A}$ is said
to be \emph{tracial} if $\varphi(ab)=\varphi(ba)$ for all $a,b\in
\mathcal{A}$.

\begin{lemma} \label{tr:lemma}
When $m$ is even, there is a unique tracial state $tr$ on
$Cl^{\Ce}(V)$. When $m$ is odd, let
$\gamma:=i^{\frac{m+1}{2}}e_1\cdots e_m$ be the \emph{chirality
operator}, where $e_1, \cdots{}, e_m$ is an orthonormal basis of
$V$. Then $\gamma$ is fixed under the action of $SO(V)$
(equivalently, $\gamma$ doesn't depend on the choice of the
ordered orthonormal basis $e_1, \cdots{}, e_m$), and there is a
unique tracial state $tr$ on $Cl^{\Ce}(V)$ such that
$tr(\gamma)=0$. In both cases, $tr$ is $SO(V)$-invariant.
\end{lemma}
\begin{proof} In both cases, the $SO(V)$-invariance of $tr$ follows from the
uniqueness. So we just need to show the uniqueness of $tr$.

When $m$ is even, $Cl^{\Ce}(V)$ is
isomorphic to the $C^*$-algebra of $2^{\frac{m}{2}}$ by $2^{\frac{m}{2}}$
matrices \cite[Theorem 1.3.2]{GM91}. The uniqueness of $tr$
%
follows from the fact that for any $n\in \Ne$ the $C^*$-algebra of $n$ by $n$ matrices has a
unique tracial state \cite[Example 8.1.2]{KR97b}.

Assume that $m$ is odd now. Then $Cl^{\Ce}(V)$ is
isomorphic to the direct sum of two copies of
the $C^*$-algebra of $2^{\frac{m-1}{2}}$ by $2^{\frac{m-1}{2}}$
matrices \cite[Theorem 1.3.2]{GM91}.
Say $Cl^{\Ce}(V)=\mathcal{A}_1\oplus \mathcal{A}_2$, where both
$\mathcal{A}_1$ and $\mathcal{A}_2$ are isomorphic to the
$C^*$-algebra of $2^{\frac{m-1}{2}}$ by $2^{\frac{m-1}{2}}$
matrices. Let $p_j$ be the projection of $Cl^{\Ce}(V)$ to
$\mathcal{A}_j$,
and let $\varphi_j$ be the unique tracial state of
$\mathcal{A}_j$. Then the tracial states of $Cl^{\Ce}(V)$ are exactly
$\lambda\varphi_1\circ p_1+ (1-\lambda)\varphi_2\circ p_2$
for $0\le \lambda\le 1$.
It is easily verified that $\gamma$ belongs to
the center of $Cl^{\Ce}(V)$. So $\gamma$ must be in
$\Ce\cdot 1_{\mathcal{A}_1}+\Ce \cdot 1_{\mathcal{A}_2}$. It's also
clear that $\gamma^2=1$ and $\gamma\not\in \Ce$. So $\gamma$ must be
$\pm (1_{\mathcal{A}_1}-1_{\mathcal{A}_2})$. It follows immediately
that $Cl^{\Ce}(V)$ has a unique tracial state $tr$ satisfying
$tr(\gamma)=0$, namely,
$tr=\frac{1}{2}(\varphi_1\circ p_1+ \varphi_2\circ p_2)$.
It is easy to check that $\gamma$ is fixed under the action
of $SO(V)$.
\end{proof}

There is a natural injective map $V\hookrightarrow  Cl(V)$.
So one may think of $V$ as a subspace of $Cl(V)$. The $C^*$-algebra norm on $Cl^{\Ce}(V)$ extends
the norm on $V$ induced from the inner product (see \cite[Theorem 1.7.22(iv)]{GM91} for the corresponding
statement for the real $C^*$-algebra norm; the proofs are similar).
 Let $M$ be an oriented Riemannian manifold of dimension $m$.
Then we have the smooth algebra bundles $ClM$ and $Cl^{\Ce}M$ over
$M$ with fibre algebras $Cl(TM_x)$ and $Cl^{\Ce}(TM_x)$
respectively, where $TM_x$ is the tangent space at $X\in M$. These
are called the \emph{Clifford algebra bundle} and the
\emph{complexified Clifford algebra bundle}. Since $TM_x\subseteq
Cl(TM_x)$, the complexified tangent bundle $TM^{\Ce}$ is a
subbundle of $Cl^{\Ce}M$. Since $Cl^{\Ce}(TM_x)$ is unital,
$C^{\infty}(M)$ is a subalgebra of $C^{\infty}(M, Cl^{\Ce}M)$.

\subsection{Compact quantum metric
spaces}\label{cqms:sub}

 Finally, we review Rieffel's theory of compact quantum
metric spaces \cite{Rieffel98b, Rieffel99b, Rieffel00, Rieffel03}.
Though Rieffel has set up his theory in the general framework of
order-unit spaces, we shall need it only for $C^*$-algebras. See
the discussion preceding Definition 2.1 in \cite{Rieffel00} for
the reason of requiring the reality condition (\ref{real:eq})
below.

\begin{definition}\cite[Definition 2.1]{Rieffel00} \label{C*QCM:def}
By a \emph{$C^*$-algebraic compact quantum metric space} we mean a
pair $(\mathcal{A}, L)$ consisting of a unital $C^*$-algebra
$\mathcal{A}$ and  a (possibly $+\infty$-valued) seminorm $L$ on
$\mathcal{A}$ satisfying
 the \emph{reality condition}
\begin{eqnarray} \label{real:eq}
L(a)&=& L(a^*)
\end{eqnarray}
for all $a\in \mathcal{A}$, such that $L$ vanishes exactly on
$\Ce$ and the metric $\rho_L$ on the state space $S(\mathcal{A})$
defined by (\ref{Lip to dist:eq}) induces the $w^*$-topology.
The \emph{radius}  of $(\mathcal{A}, L)$ is defined to be the radius of $(S(\mathcal{A}), \rho_L)$.
We say that $L$ is a \emph{Lip-norm}.
\end{definition}

Let $\mathcal{A}$ be a unital $C^*$-algebra and let $L$ be a
(possibly $+\infty$-valued) seminorm on $\mathcal{A}$ vanishing on
$\Ce$. Then $L$ and $\pa \cdot \pa$ induce (semi)norms ${\tilde
L}$ and $\pa \cdot\pa^{\sim}$ respectively on the quotient space
$\tilde{\mathcal{A}}=\mathcal{A}/\Ce$.

\begin{notation} \label{ball:notation}
For any $r\ge 0$, let
\begin{eqnarray*}
\mathcal{D}_{r}(\mathcal{A}):=\{a\in \mathcal{A}: L(a)\le 1, \pa a\pa \le r\}.
\end{eqnarray*}
\end{notation}

The main criterion for when a seminorm $L$ is a Lip-norm is the following:

\begin{proposition}\cite[Proposition 1.6, Theorem 1.9]{Rieffel98b}
\label{criterion of Lip:prop}
Let $\mathcal{A}$ be a unital $C^*$-algebra and let
$L$ be a (possibly $+\infty$-valued) seminorm on $\mathcal{A}$
satisfying the reality condition (\ref{real:eq}).
Assume that $L$ takes finite values on a dense
subspace of $\mathcal{A}$, and that $L$ vanishes exactly on $\Ce$.
Then $L$ is a Lip-norm if and only if

%
\noindent \, \, \, \, \,\, \, \,
(1) there is a constant $K\ge 0$ such that $\pa \cdot\pa^{\sim}\le K
\tilde{L}$ on $\tilde{\mathcal{A}}$;

\noindent
and \, \, \, \,
(2) for any $r\ge 0$, the  ball $\mathcal{D}_{r}(\mathcal{A})$
is totally bounded in $\mathcal{A}$ for $\pa \cdot \pa$;

\noindent \, \, \, \, \,
or
(2') for some $r> 0$, the  ball $\mathcal{D}_{r}(\mathcal{A})$
is totally bounded in $\mathcal{A}$ for $\pa \cdot \pa$.

In this event, $r_{\mathcal{A}}$ is exactly the minimal $K$ such
that $\pa \cdot\pa^{\sim}\le K \tilde{L}$ on
$(\tilde{\mathcal{A}})_{sa}$.
\end{proposition}

\section{Connes and Dubois-Violette's formulation of $\theta$-deformations}
\label{ConsConnesDubois:sec}

Though the Dirac operator does not depend on $\theta$ in Connes
and Landi's formulation of $\theta$-deformations in \cite[Section
5]{CL01}, it does in Connes and Dubois-Violette's formulation in
\cite{CD01}. In this section we review the formulation of
$\theta$-deformations by Connes and Dubois-Violette
\cite[Sections 11 and 13]{CD01}, including the deformation of both
the algebra and the Dirac operator.

Let $M$ be a smooth manifold with a smooth action $\sigma_M$ of
$\mathbb{T}^n$. We denote by  $\sigma$ the induced action of
$\mathbb{T}^n$ on the {\LCSA} $C^{\infty}(M)$. Then $\sigma$ is
continuous. By Proposition~\ref{tensor prod of alg:lemma} the
tensor product completion $C^{\infty}(M)\hat{\otimes}
\mathcal{A}^{\infty}_{\theta}$ is a {\LCSA}.
The tensor product action $\sigma \hat{\otimes} \tau^{-1}$ of
$\mathbb{T}^n$ on $C^{\infty}(M)\hat{\otimes}
\mathcal{A}^{\infty}_{\theta}$ is also continuous. The deformed
smooth algebra \cite[Section 11]{CD01}, denoted by
$C^{\infty}(M_{\theta})$, is then defined as the fixed-point space
of this action, \ie
$C^{\infty}(M_{\theta})=(C^{\infty}(M)\hat{\otimes}
\mathcal{A}^{\infty}_{\theta})^{\sigma \hat{\otimes} \tau^{-1}}$.
Clearly, this is a {\LCSA}.

Suppose $M$ is equipped with a $\sigma_M$-invariant Riemannian
metric. (For any Riemannian metric on $M$, we can always
  integrate it over $\mathbb{T}^n$ to make it
  $\sigma_M$-invariant.)
Also assume that $M$ is a spin manifold and that $\sigma_M$ lifts to a
smooth action
$\sigma_S$ of $\mathbb{T}^n$ on the spin bundle $S$, \ie the following
diagram
\begin{equation*}
\begin{CD}
 S  @>>\sigma_{S, x}>  S \\
@VVV     @VVV \\
 M  @>>\sigma_{M, x}> M
\end{CD}
\end{equation*}
is commutative for every $x\in \mathbb{T}^n$.
(Usually $\sigma_M$ doesn't lift directly to $S$, but lifts
  only modulo $\pm I$, \ie there is a twofold covering
  $\mathbb{T}^n\rightarrow \mathbb{T}^n$ such that $\sigma_M$ lifts to
  an action of the two-folding covering on $S$. Correspondingly, Connes
  and Dubois-Violette defined the various deformed structures using
  tensor product with $\mathcal{A}_{\frac{1}{2}\theta}$ instead of
  $\mathcal{A}_{\theta}$. But for the
  deformed algebras and Dirac operators, the
  difference is just a matter of parameterization.)
We denote the induced continuous action of $\mathbb{T}^n$ on $C^{\infty}(M, S)$
also by $\sigma$. Then $C^{\infty}(M, S)$ is a locally convex
left $C^{\infty}(M)$-module
and
\begin{eqnarray*}
\sigma_x(f\psi)=\sigma_x(f)\sigma_x(\psi)
\end{eqnarray*}
for all $f\in C^{\infty}(M), \psi \in C^{\infty}(M, S)$ and $x\in
\mathbb{T}^n$. We also have the tensor product completion
$C^{\infty}(M, S)\hat{\otimes} \mathcal{A}^{\infty}_{\theta}$,
which is a locally convex left module over
$C^{\infty}(M)\hat{\otimes} \mathcal{A}^{\infty}_{\theta}$ by
Proposition~\ref{tensor prod of alg:lemma}. The tensor product
action $\sigma \hat{\otimes} \tau^{-1}$ of $\mathbb{T}^n$ on
$C^{\infty}(M, S)\hat{\otimes} \mathcal{A}^{\infty}_{\theta}$ is
still continuous. The deformed spin bundle, denoted by
$C^{\infty}(M_{\theta}, S)$, is then defined as the fixed-point
space of this action, \ie $C^{\infty}(M_{\theta},
S)=(C^{\infty}(M, S)\hat{\otimes}
\mathcal{A}^{\infty}_{\theta})^{\sigma \hat{\otimes} \tau^{-1}}$.
This is a locally convex left $C^{\infty}(M_{\theta})$-module. Let
$D$ be the Dirac operator on $C^{\infty}(M, S)$. This is a
first-order linear differential operator. So it is easy to see
that $D$ is continuous with respect to the locally convex topology
on $C^{\infty}(M, S)$. Then we have the tensor product linear map
$D\hat{\otimes} I$ from $C^{\infty}(M, S)\hat{\otimes}
\mathcal{A}^{\infty}_{\theta}$ to itself.
 Notice that $D$ commutes with the action $\sigma$, so
$D\hat{\otimes} I$ commutes with the action $\sigma \hat{\otimes}
\tau^{-1}$. Therefore $C^{\infty}(M_{\theta}, S)$ is stable under
$D\hat{\otimes} I$.
Denote by $D_{\theta}$ the restriction of $D\hat{\otimes} I$ to
$C^{\infty}(M_{\theta}, S)$.

Assume further that $M$ is compact. As usual,
one defines a positive-definite scalar product on $C^{\infty}(M, S)$ by
\begin{eqnarray*}
<\psi, \psi'>=\int_M(\psi, \psi')\, vol,
\end{eqnarray*}
where $vol$ is the Riemannian volume form. Denote by
$\mathcal{H}=L^2(M, S)$ the Hilbert space obtained by completion.
Then $C(M)$ has a natural faithful representation on $\mathcal{H}$
by multiplication, and we shall think of $C(M)$ as a subalgebra of
$B(\mathcal{H})$, the $C^*$-algebra of all bounded operators on
$\mathcal{H}$. The action $\sigma$ uniquely extends to a
continuous unitary representation of $\mathbb{T}^n$ in
$\mathcal{H}$, which will be still denoted by $\sigma$. On the
other hand, $\mathcal{A}_{\theta}$ has an inner product induced by
the unique $\tau$-invariant tracial state. Denote by
$L^2(\mathcal{A}_{\theta})$ the Hilbert space obtained by
completion. Then $\mathcal{A}_{\theta}$ acts on
$L^2(\mathcal{A}_{\theta})$ faithfully by the $GNS$ construction,
and we shall also think of $\mathcal{A}_{\theta}$ as a subalgebra
of $B(L^2(\mathcal{A}_{\theta}))$. The action $\tau$ also extends
to a continuous unitary representation of $\mathbb{T}^n$ in
$L^2(\mathcal{A}_{\theta})$. Let $\mathcal{H}\bar{\otimes}
L^2(\mathcal{A}_{\theta})$ be the Hilbert space tensor product.
Then we have the continuous tensor product action $\sigma
\bar{\otimes} \tau^{-1}$ on $\mathcal{H}\bar{\otimes}
L^2(\mathcal{A}_{\theta})$. The deformed Hilbert space, denoted by
$\mathcal{H}_{\theta}$, is defined as the fixed-point space of
$\mathcal{H}\bar{\otimes} L^2(\mathcal{A}_{\theta})$ under the
action $\sigma \bar{\otimes} \tau^{-1}$. Clearly the maps
$C^{\infty}(M, S)\rightarrow \mathcal{H}$ and
$\mathcal{A}^{\infty}_{\theta}\rightarrow
L^2(\mathcal{A}_{\theta})$ are continuous with respect to the
locally convex topologies on $C^{\infty}(M, S),\,
\mathcal{A}^{\infty}_{\theta}$ and the norm topologies on $
\mathcal{H},\, L^2(\mathcal{A}_{\theta})$. Then we have the
sequence of continuous linear maps
\begin{eqnarray*}
C^{\infty}(M,
S)\hat{\otimes}\mathcal{A}^{\infty}_{\theta}\overset{\phi}\rightarrow
\mathcal{H}\hat{\otimes}_{\pi}L^2(\mathcal{A}_{\theta})\overset{\psi}\rightarrow
\mathcal{H}\bar{\otimes} L^2(\mathcal{A}_{\theta}),
\end{eqnarray*}
where $\mathcal{H}\hat{\otimes}_{\pi}L^2(\mathcal{A}_{\theta})$ is
the completion of the projective tensor product of $\mathcal{H}$
and $L^2(\mathcal{A}_{\theta})$. Let $\Phi:C^{\infty}(M,
S)\hat{\otimes}\mathcal{A}^{\infty}_{\theta}\rightarrow
\mathcal{H}\bar{\otimes} L^2(\mathcal{A}_{\theta})$ be the
composition. Then $\Phi$ is $\mathbb{T}^n$-equivariant. So $\Phi$
maps $C^{\infty}(M_{\theta}, S)$ into $\mathcal{H}_{\theta}$. Let
$\Phi_{\theta}$ be the restriction of $\Phi$ to
$C^{\infty}(M_{\theta}, S)$.

\begin{lemma} \label{injective Phi:lemma}
Both maps $\phi:C^{\infty}(M,
S)\hat{\otimes}\mathcal{A}^{\infty}_{\theta}\rightarrow
\mathcal{H}\hat{\otimes}_{\pi}L^2(\mathcal{A}_{\theta})$ and
$\psi: \\
\mathcal{H}\hat{\otimes}_{\pi}L^2(\mathcal{A}_{\theta})\rightarrow
\mathcal{H}\bar{\otimes} L^2(\mathcal{A}_{\theta})$ are injective.
Consequently, $\Phi$ and $\Phi_{\theta}$ are injective.
\end{lemma}
\begin{proof}
We'll prove the injectivity of $\phi$. The proof for $\psi$ is
similar. Recall the notation at the end of
Section~\ref{intro:sec}.
We shall need the following well-known fact
several times. We omit the proof.
\begin{lemma} \label{proj TVS:lemma}
Let $G$ be a compact group.
Let $\alpha$ be a continuous action of $G$ on a complex
complete {\LCTVS}
 $V$.  For a continuous
$\Ce$-valued function $\varphi$ on $G$ let
\begin{eqnarray*}
\alpha_{\varphi}(v)=\int_G \varphi(x)\alpha_x(v)\, dx
\end{eqnarray*}
for $v\in V$. Then $\alpha_{\varphi}:V\rightarrow V$ is a continuous linear
map. If $\mathcal{J}$ is a finite subset of $\hat{G}$ and if $\varphi$ is
a linear combination of the characters of $\gamma\in
\bar{\mathcal{J}}$, then $\alpha_{\varphi}(V)\subseteq V_{\mathcal{J}}$.
Let
\begin{eqnarray*}
\alpha_{\mathcal{J}}=\alpha_{\sum_{\gamma \in \mathcal{J}}\dim(\gamma)
\overline{\chi_{\gamma}}}.
\end{eqnarray*}
(When $\mathcal{J}$ is a one-element set $\{ \gamma\}$, we'll simply
write $\alpha_{\gamma}$ for $\alpha_{\{\gamma\}}$.)
Then
$\alpha_{\mathcal{J}}(v)=v$ for all $v\in V_{\mathcal{J}}$, and
$\alpha_{\mathcal{J}}(v)=0$ for all $v\in V_{\gamma}$ with $\gamma \in
\hat{G}\setminus \mathcal{J}$.
\end{lemma}

From Proposition~\ref{linear of integ:prop} we also have:
\begin{lemma} \label{exchange integ:lemma}
Let $G$ be a compact group with continuous actions $\alpha$ and
$\beta$ on complex complete {\LCTVS}
 $V$ and $W$.
Let $\phi:V\rightarrow W$ be a continuous
 $G$-equivariant linear map, and let $\varphi:G\rightarrow \Ce$ be a
 continuous function.
Then
\begin{eqnarray} \label{exchange integ:eq}
\phi\circ \alpha_{\varphi}=\beta_{\varphi}\circ \phi.
\end{eqnarray}
In particular, let $\mathcal{J}$ be a finite subset of $\hat{G}$. Then
\begin{eqnarray*}
\phi\circ \alpha_{\mathcal{J}}=\beta_{\mathcal{J}}\circ \phi.
\end{eqnarray*}
\end{lemma}

We shall need the following lemma a few times:
\begin{lemma} \label{tensor limit vanish:lemma}
Let $G$ be a compact group, and let $h$ be a continuous $\Ce$-valued function on $G$ with $h(e_G)=0$.
Then for any $\epsilon>0$
there is a nonnegative function $\varphi$ on $G$ such that
$\varphi$ is a linear combination of finitely many characters,
$\pa \varphi\pa_1=1$, and $\pa \varphi\cdot h\pa_1< \epsilon$.
\end{lemma}
\begin{proof} Notice that the left regular representation of
$G$ on $L^2(G)$ is faithful. Since the left regular representation
is a Hilbert space direct sum of irreducible representations, we
see that any $x\neq e_G$ acts nontrivially in some $\gamma \in
\hat{G}$. Let $\mathcal{U}$ be an open neighborhood of $e_G$ such
that $|h(x)|< \epsilon/2$ for all $x\in \mathcal{U}$. For any
$x\in G\setminus \mathcal{U}$, suppose that $x$ acts nontrivially
in $\gamma_x\in \hat{G}$. Then there is some open neighborhood
$\mathcal{U}_x$ of $x$ such that $x'$ acts nontrivially in
$\gamma_x$ for all $x'\in \mathcal{U}_X$. Since $G\setminus
\mathcal{U}$ is compact, we can find $x_1, {\cdots}, x_m \in
G\setminus \mathcal{U}$ so that $U_{x_1}, {\cdots}, U_{x_m}$ cover
$G\setminus \mathcal{U}$. Let
$\mathcal{J}_{\mathcal{U}}=\{\gamma_{x_1}, {\cdots},
\gamma_{x_m}\}$. Then no element in $G\setminus \mathcal{U}$ acts
trivially in all $\gamma \in \mathcal{J}_{\mathcal{U}}$. Let
$\pi_1$ be the direct sum of one copy for each $\gamma$ in
$\mathcal{J}_{\mathcal{U}}\cup \{\gamma_0\}$, and let
$\chi_{\pi_1}$ be the character of $\pi_1$.

Let $\pi=\pi_1\otimes \overline{\pi_1}$.
Also let  $\chi$ be the character of $\pi$.
Note that $\chi(x)=|\chi_{\pi_1}(x)|^2\ge 0$ for all $x\in G$.
Let $\varphi_n=\chi^n/\pa \chi^n\pa_1$. Then each $\varphi_n$ is
a linear combination of finitely many characters.
Since every element in $G\setminus \mathcal{U}$ acts nontrivially in $\pi$,
$\chi(x)<\chi(e_G)$ on $G\setminus \mathcal{U}$.
Therefore it's easy to see (cf. the proof of Theorem 8.2 in
\cite{Rieffel00}) that
$\int_{G\setminus \mathcal{U}}\varphi_n(x)\, dx\to 0$ as $n\to \infty$,
and hence
\begin{eqnarray*}
\limsup_{n\to +\infty}\int_G|\varphi_n(x)h(x)|\, dx \le
 \sup_{x\in \mathcal{U}}|h(x)|< \epsilon.
\end{eqnarray*}
So when $n$ is big enough, we have that
$\pa \varphi_n\cdot h\pa_1<\epsilon$.
\end{proof}

As a corollary of Lemma~\ref{tensor limit vanish:lemma} we have:
\begin{lemma} \label{comp vanish to vanish:lemma}
Let $G$ be a compact group.
Let $\alpha$ be a continuous action of $G$ on
a complex complete {\LCTVS} $V$. Let $v\in V$. If
$\alpha_{\gamma}(v)=0$ for all $\gamma \in \hat{G}$, then $v=0$.
\end{lemma}
\begin{proof}
Let $\mathfrak{p}$ be a continuous seminorm on $V$, and let
$\epsilon>0$.
Define a function $h$ on $G$ by
 $h(x)=\mathfrak{p}(v-\alpha_x(v))$. Then
$h$ is continuous on $G$, and $h(e_G)=0$.
Pick $\varphi$ for $h$ and $\epsilon$ in
Lemma~\ref{tensor limit vanish:lemma}. According to the assumption we
have $\alpha_{\varphi}(v)=0$. Then
\begin{eqnarray*}
\mathfrak{p}(v)=\mathfrak{p}(v-\alpha_{\varphi}(v))=
\mathfrak{p}(\int_G\varphi(x)(v-\alpha_x(v))\, dx)
\le \int_G\varphi(x)
h(x)\, dx <\epsilon.
\end{eqnarray*}
Since the topology on $V$ is defined by all the continuous seminorms,
we see that $v=0$.
\end{proof}

We are ready to prove Lemma~\ref{injective Phi:lemma}. Let
$\alpha=I\hat{\otimes}\tau$ acting on $V=C^{\infty}(M,
S)\hat{\otimes}\mathcal{A}^{\infty}_{\theta}$, and let
$\beta=I\hat{\otimes}\tau$ acting on
$\mathcal{H}\hat{\otimes}_{\pi}L^2(\mathcal{A}_{\theta})$. Let
$\phi$ be as in Lemma~\ref{injective Phi:lemma}. Then $\phi\circ
\alpha=\beta\circ \phi$. Recall the notation about
$\mathcal{A}_{\theta}$ in subsection~\ref{LCTVS:sub}. For any
$q\in \Ze^n=\widehat{\mathbb{T}^n}$ clearly $\alpha_q$ maps
$C^{\infty}(M)\otimes_{alg} \mathcal{A}^{\infty}_{\theta}$ onto
$C^{\infty}(M)\otimes u_q$. Since $\alpha_q$ is continuous, by
Lemma~\ref{proj TVS:lemma} it follows immediately that
$V_q=\alpha_q(V)= C^{\infty}(M, S)\otimes u_q$. Let $f\in
ker(\phi)$. For any $q\in \Ze^n$ by Lemma~\ref{exchange
integ:lemma}
 $\phi(\alpha_q(f))=\beta_q(\phi(f))=0$.
Now $\alpha_q(f)\in C^{\infty}(M, S)\otimes u_q$, and clearly $\phi$
restricted to $C^{\infty}(M, S)\otimes u_q$ is injective.
So $\alpha_q(f)=0$.
From Lemma~\ref{comp vanish to vanish:lemma} we see that $f=0$.
\end{proof}

\begin{lemma} \label{dense Hilbert:lemma}
The image $\Phi_{\theta}(C^{\infty}(M_{\theta}, S))$ is dense in
$\mathcal{H}_{\theta}$.
\end{lemma}

Clearly $\Phi(C^{\infty}(M,
S)\hat{\otimes}\mathcal{A}^{\infty}_{\theta})$ is dense in
$\mathcal{H}\bar{\otimes} L^2(\mathcal{A}_{\theta})$, so this is
an immediate consequence of the following:
\begin{lemma} \label{dense:lemma}
Let $G$ be a compact group.
Let  $\alpha$ and $\beta$ be continuous actions of $G$ on
complex complete {\LCTVS}
 $V$ and $W$ respectively. Let $\phi:V\rightarrow W$ be a continuous
 $G$-equivariant linear map such that $\phi(V)$ is dense in $W$.
Then $\phi(V^{\alpha})$ is dense in $W^{\beta}$.
\end{lemma}
\begin{proof}
Recall that $\gamma_0$ is the trivial representation of $G$.
By Lemma~\ref{proj TVS:lemma} $\beta_{\gamma_0}$ is continuous.
So $\beta_{\gamma_0}(\phi(V))$ is dense in
$\beta_{\gamma_0}(W)=W^{\beta}$. But
$\beta_{\gamma_0}(\phi(V))=\phi(\alpha_{\gamma_0}(V))=\phi(V^{\alpha})$
according to Lemma~\ref{exchange integ:lemma}. The conclusion follows.
\end{proof}

The Dirac operator $D$ is essentially self-adjoint on $\mathcal{H}$
\cite[Theorem 5.7]{LM89}. Then
$D\otimes I$ is also essentially self-adjoint on
$\mathcal{H}\bar{\otimes} L^2(\mathcal{A}_{\theta})$
\cite[Proposition 11.2.37]{KR97b}. Denote
 its closure by $D^{L^2}$.
\begin{lemma} \label{Phi in D:lemma}
$\Phi(C^{\infty}(M, S)\hat{\otimes}\mathcal{A}^{\infty}_{\theta})$
is contained in the domain of $D^{L^2}$, and
\begin{eqnarray} \label{D Psi:eq}
D^{L^2}\circ \Phi=\Phi\circ (D\hat{\otimes}I).
\end{eqnarray}
\end{lemma}
\begin{proof}
For any $y\in C^{\infty}(M,
S)\hat{\otimes}\mathcal{A}^{\infty}_{\theta}$, take a net $y_j$ in
$C^{\infty}(M, S)\otimes_{alg}\mathcal{A}^{\infty}_{\theta}$
converging to $y$. Then $\Phi(y_j)\to \Phi(y), \,
(D\hat{\otimes}I)(y_j)\to (D\hat{\otimes}I)(y)$ and
$D^{L^2}(\Phi(y_j))=\Phi((D\hat{\otimes}I)(y_j))\to
\Phi((D\hat{\otimes}I)(y))$. So $\Phi(y)$ is contained in the
domain of $D^{L^2}$, and
$D^{L^2}(\Phi(y))=\Phi((D\hat{\otimes}I)(y))$.
\end{proof}
So the intersection of $\mathcal{H}_{\theta}$ and the domain of
$D^{L^2}$ contains $\Phi_{\theta}(C^{\infty}(M_{\theta}, S))$,
which is dense in $\mathcal{H}_{\theta}$ by Lemma~\ref{dense
Hilbert:lemma}. Clearly $D\otimes I$ commutes with the action
$\sigma\bar{\otimes}\tau^{-1}$, and thus so does $D^{L^2}$. Hence
$D^{L^2}$ maps the intersection of $\mathcal{H}_{\theta}$ and the
domain of $D^{L^2}$ into $\mathcal{H}_{\theta}$.
Therefore the
restriction of $D^{L^2}$ to $\mathcal{H}_{\theta}$ is also
self-adjoint. The deformed Dirac operator, denoted by
$D^{L^2}_{\theta}$, is then defined to be this restriction.

Similarly, the maps $C^{\infty}(M)\rightarrow C(M)$ and
$\mathcal{A}^{\infty}_{\theta}\rightarrow \mathcal{A}_{\theta}$
are continuous with respect to the locally convex topologies on
$C^{\infty}(M),\, \mathcal{A}^{\infty}_{\theta}$ and the norm
topologies on $C(M),\, \mathcal{A}_{\theta}$. So we have the
$\mathbb{T}^n$-equivariant continuous linear map
$\Psi:C^{\infty}(M)\hat{\otimes}\mathcal{A}^{\infty}_{\theta}
\rightarrow C(M)\otimes \mathcal{A}_{\theta}$ , where $C(M)\otimes
\mathcal{A}_{\theta}$ is the spatial $C^*$-algebraic tensor
product of $C(M)$ and $\mathcal{A}_{\theta}$ \cite[Appendix
T.5]{Wegge93}.

\begin{definition} \label{deformed cont:def}
We define the deformed continuous algebra, $C(M_{\theta})$, to be
the fixed-point algebra
$(C(M)\otimes \mathcal{A}_{\theta})^{\sigma \otimes \tau^{-1}}$.
\end{definition}
Then $\Psi$ maps $C^{\infty}(M_{\theta})$
into $C(M_{\theta})$. By similar arguments as
in Lemma~\ref{injective Phi:lemma} and \ref{dense:lemma} we have
\begin{lemma} \label{injective Psi:lemma}
The map $\Psi$ is injective, and
$\Psi(C^{\infty}(M_{\theta}))$ is dense in $C(M_{\theta})$.
\end{lemma}

Clearly $\mathcal{H}_{\theta}$ is stable under the action of
elements in $C(M_{\theta})$. So we can define
$\Psi_{\theta}:C^{\infty}(M_{\theta})\rightarrow
B(\mathcal{H}_{\theta})$ as the composition of
$C^{\infty}(M_{\theta}) \rightarrow C(M_{\theta})$ and the
restriction map of $C(M_{\theta})$ to $B(\mathcal{H}_{\theta})$.
We shall see later in Proposition~\ref{same norm theta:prop} that
the restriction map of $C(M_{\theta})$ to
$B(\mathcal{H}_{\theta})$ is isometric. So we may also think of
$C(M_{\theta})$ as a subalgebra of $B(\mathcal{H}_{\theta})$.
 Then the closure of
$\Psi_{\theta}(C^{\infty}(M_{\theta}))$ is just
$C(M_{\theta})$.

We shall see later in Proposition~\ref{commutator psi:prop}
that  the domain of $D^{L^2}_{\theta}$ is stable under
$\Psi_{\theta}(f)$, and that the commutator $[D^{L^2}_{\theta},
\Psi_{\theta}(f)]$
is bounded for every $f\in C^{\infty}(M_{\theta})$.

\begin{definition} \label{L_{theta}:def}
We define the \emph{deformed Lip-norm}, denoted by
$L_{\theta}$, on $C(M_{\theta})$ by
\begin{eqnarray*}
L_{\theta}(f):=\begin{cases}
\pa [D^{L^2}_{\theta}, f]\pa, & \mbox{ if } f\in
\Psi_{\theta}(C^{\infty}(M_{\theta}));\\
+\infty, & \mbox{ otherwise }.
\end{cases}
\end{eqnarray*}
\end{definition}


\section{Lip-norms and Compact Group Actions}
\label{Lip action:sec}

In this section we consider a general situation in which there are a
seminorm and a compact group action. We show that under certain
compatibility hypotheses we can use this group action to prove that
the seminorm is a Lip-norm. The strategy is a generalization of
the one Rieffel used to deal
with Lip-norms associated to ergodic compact (Lie) group actions
\cite{Rieffel98b, Rieffel00}. We'll
see that $\theta$-deformations
fit into this general picture.

Throughout this section we assume that $G$ is an arbitrary
compact group which has a fixed length function
$\mathnormal{l}$,  \ie a continuous
real-valued function, $\mathnormal{l}$, on $G$ such that
\begin{eqnarray*}
\mathnormal{l}(xy)&\le &\mathnormal{l}(x)+\mathnormal{l}(y)\mbox{ for
  all } x, y\in G \\
\mathnormal{l}(x^{-1})&=& \mathnormal{l}(x) \mbox{ for all }x\in G \\
\mathnormal{l}(x)&=& 0 \mbox{ if and only if } x=e_G,
\end{eqnarray*}
where $e_G$ is the identity of $G$.

\begin{theorem} \label{use of action:thm}
Let $\mathcal{A}$ be a unital $C^*$-algebra, let
$L$ be a (possibly $+\infty$-valued) seminorm on $\mathcal{A}$
satisfying the reality condition (\ref{real:eq}),
and let $\alpha$ be a strongly continuous action of $G$
on $\mathcal{A}$. Assume that $L$ takes finite values on a dense
subspace of $\mathcal{A}$, and that $L$ vanishes on $\Ce$.
Let $L^{\mathnormal{l}}$ be the (possibly $+\infty$-valued)
seminorm on $\mathcal{A}$ defined by
\begin{eqnarray} \label{def of L:eq}
L^{\mathnormal{l}}(a)=\sup \{\frac{\pa \alpha_x(a)-a\pa}{\mathnormal{l}(x)}: x\in G, x\neq e_G\}.
\end{eqnarray}
Suppose that the following conditions are
satisfied:

(1) there is some constant $C>0$ such that $L^{\mathnormal{l}}\le C\cdot L$ on
    $\mathcal{A}$;

(2) for any linear combination $\varphi$ of finitely many characters
    on $G$ we have $L\circ \alpha_{\varphi}\le$ $\pa \varphi\pa_1\cdot
    L$ on $\mathcal{A}$, where
    $\alpha_{\varphi}$ is the linear map on
    $\mathcal{A}$ defined in Lemma~\ref{proj TVS:lemma};

(3) for each $\gamma\in \hat{G}$ with $\gamma\neq \gamma_0$ the ball
    $\mathcal{D}_{r}(\mathcal{A}_{\gamma}):=\{a\in
    \mathcal{A}_{\gamma}:L(a)\le 1, \pa a\pa \le r\}$ is totally
    bounded for some $r>0$, and the only element in $A_{\gamma}$
    vanishing under $L$ is $0$;

(4) there is a unital $C^*$-algebra $\mathcal{B}$  containing
$\mathcal{A}_{\gamma_0}=\mathcal{A}^{\alpha}$, with
a Lip-norm $L_{\mathcal{B}}$, such
  that $L_{\mathcal{B}}$ extends the restriction of $L$ to $\mathcal{A}_{\gamma_0}$.

Then $(\mathcal{A}, L)$ is a $C^*$-algebraic compact quantum
metric space with $r_{\mathcal{A}}\le r_{\mathcal{B}}+
C\int_Gl(x)\, dx$.
\end{theorem}

\begin{remark} \label{use of action:remark}
(1) We assume the existence of $(\mathcal{B}, L_{\mathcal{B}})$ in the
    condition (4) only for the convenience of application. In fact,
    conditions (2) and (4) imply that $L$ restricted to
    $\mathcal{A}_{\gamma_0}$ is a Lip-norm on $\mathcal{A}_{\gamma_0}$:
    for any $a\in \mathcal{A}_{\gamma_0}$ and $\epsilon>0$ pick $a'\in
    \mathcal{A}$ with $L(a')<\infty$ and $\pa a-a'\pa<\epsilon$.
    Then by Lemma~\ref{proj TVS:lemma}
    $\alpha_{\gamma_0}(a')\in \mathcal{A}_{\gamma_0}$ and $\pa
    a-\alpha_{\gamma_0}(a')\pa=\pa
    \alpha_{\gamma_0}(a-a')\pa<\epsilon$.
    By the condition (2) $L(\alpha_{\gamma_0}(a'))<\infty$.
    Therefore $L$ takes finite values on a dense subspace of
    $\mathcal{A}_{\gamma_0}$.
    Then from Proposition~\ref{criterion of Lip:prop} it is easy
    to see that $L$ restricted to $\mathcal{A}_{\gamma_0}$
     is a Lip-norm on $\mathcal{A}_{\gamma_0}$. Consequently, we
    may take $\mathcal{B}$ to be $\mathcal{A}_{\gamma_0}$ itself.

(2) Conditions (1) and (2) in Theorem~\ref{use of action:thm} enable us to
reduce the study of $L$ to that of the restriction of $L$ to each
$\mathcal{A}_{\gamma}$. Conditions (3) and (4) say roughly that $L$
restricted to each $\mathcal{A}_{\gamma}$ is a Lip-norm.

(3) Usually it is not hard to verify the condition (2).
In particular, by Lemma~\ref{lower (2):lemma} it holds
when $L$ is $\alpha$-invariant and
lower semicontinuous on $\{a\in \mathcal{A}:L(a)<+\infty\}$, and
$\{a\in \mathcal{A}:L(a)<+\infty\}$ is stable under $\alpha_{\gamma}$
for every $\gamma \in \hat{G}$.
\end{remark}

\begin{lemma} \label{lower (2):lemma}
Let $\alpha$ be a strongly continuous action of $G$ on a $C^*$-algebra
$\mathcal{A}$, and  let $L$ be a (possibly $+\infty$-valued) seminorm on
$\mathcal{A}$. Suppose that $L$ is $\alpha$-invariant and
lower semicontinuous on $\{a\in \mathcal{A}:L(a)<+\infty\}$. For
any continuous function $\varphi:G\to \Ce$, if
$\{a\in \mathcal{A}:L(a)<+\infty\}$ is stable under
 the map $\alpha_{\varphi}:\mathcal{A}\rightarrow \mathcal{A}$ defined in
 Lemma~\ref{proj TVS:lemma}, then
\begin{eqnarray*}
L\circ \alpha_{\varphi}\le \pa \varphi\pa_1\cdot L
\end{eqnarray*}
on $\mathcal{A}$.
\end{lemma}
\begin{proof}
We only need to show $L(\alpha_{\varphi}(a))\le \pa
\varphi\pa_1\cdot L(a)$ for each $a\in \mathcal{A}$ with
$L(a)<+\infty$. But
\begin{eqnarray*}
\alpha_{\varphi}(a)=\lim_{\Delta\to 0}
\sum^k_{j=1}\alpha_{g_j}(a)\mu(E_j)\varphi(g_j),
\end{eqnarray*}
where $\mu$ is the normalized Haar measure on $G$, $(E_1, \cdots,
E_k)$ is a partition of $G$, $g_j\in E_j$,
$\Delta(E_j):=\sup\{\max(|\varphi(x)-\varphi(y)|, \,
|\alpha_x(a)-\alpha_y(a)|):x,y \in E_j\}$ and $\Delta=\max_{1\le
j\le k}\Delta(E_j)$. By the assumptions we have
\begin{eqnarray*}
L(\alpha_{\varphi}(a))&\le &\liminf_{\Delta\to 0}
L(\sum^k_{j=1}\alpha_{g_j}(a)\mu(E_j)\varphi(g_j)) \\
&\le &
L(a)  \liminf_{\Delta\to 0}
\sum^k_{j=1}\mu(E_j)|\varphi(g_j)|=L(a)\pa \varphi\pa_1.
\end{eqnarray*}
\end{proof}

For $\theta$-deformations of course  $\mathcal{A}$ is
$C(M_{\theta})$.
Notice that $\mathbb{T}^n$ has a natural action $I\otimes \tau$ on
$C(M_{\theta})$.
They will be our $G$ and $\alpha$.

The following lemma is
 a generalization of Lemmas 8.3 and 8.4 in \cite{Rieffel00}.

\begin{lemma} \label{finite approx*:lemma}
For any $\epsilon>0$ there is a finite subset
$\mathcal{J}=\bar{\mathcal{J}}$ in
$\hat{G}$, containing $\gamma_0$, depending only
on $\mathit{l}$ and $\epsilon/C$, such that for any strongly continuous
isometric action $\alpha$ on a complex Banach space $V$ with a
(possibly $+\infty$-valued) seminorm $L$ on $V$ satisfying conditions (1) and (2) (with
$\mathcal{A}$ replaced by $V$) in
Theorem~\ref{use of action:thm},
and for any $v\in V$, there is some $v'\in V_{\mathcal{J}}$ with
\begin{eqnarray*}
\pa v'\pa\le \pa v\pa, \quad L(v')\le
L(v), \quad \mbox{ and }\pa v-v'\pa\le \epsilon L(v).
\end{eqnarray*}
If $V$ has an isometric involution $*$ invariant under $\alpha$,
then when $v$ is self-adjoint
we can choose $v'$ also to be self-adjoint.
\end{lemma}
\begin{proof}
Pick $\varphi$ for $\mathnormal{l}$ and $\epsilon/C$ as in
Lemma~\ref{tensor limit vanish:lemma}. Then there is a
finite subset $\mathcal{J}\subseteq \hat{G}$ such that $\varphi$ is
a linear combination of characters $\chi_{\gamma}$ for $\gamma\in
\mathcal{J}$. Replacing $\mathcal{J}$ by $\mathcal{J}\cup
\bar{\mathcal{J}}$, we may assume that $\mathcal{J}=
\bar{\mathcal{J}}$.
For any $v\in V$ clearly
\begin{eqnarray*}
\pa \alpha_{\varphi}(v)\pa \le \pa \varphi\pa_1\cdot \pa v\pa =\pa v\pa.
\end{eqnarray*}
A simple calculation as in the proof of \cite[Lemma
8.3]{Rieffel00} tells us that
\begin{eqnarray*}
\pa v-\alpha_{\varphi}(v)\pa
&\le &L^{\mathnormal{l}}(v)\int_G\varphi(x)\mathnormal{l}(x)\, dx\le
\frac{\epsilon}{C} L^{\mathnormal{l}}(v).
\end{eqnarray*}
Then it follows from the condition (1) in Theorem~\ref{use of
  action:thm} that $\pa v-\alpha_{\varphi}(v)\pa\le \epsilon L(v)$.
Also from the condition (2) we see that $L(\alpha_{\varphi}(v))\le
  L(v)$.
So for any $v\in \mathcal{A}$, the element
$v'=\alpha_{\varphi}(v)$  satisfies the requirement.

Notice that $\varphi$ is real-valued, so when $v$ is self-adjoint, so
is $\alpha_{\varphi}(v)$.
\end{proof}

\begin{proof}[Proof of Theorem~\ref{use of action:thm}]
We verify the conditions in Proposition~\ref{criterion of
Lip:prop} for $(\mathcal{A}, L)$ to be a compact quantum metric
space one by one.

\begin{lemma} \label{condition 1*:lemma}
For any $a\in \mathcal{A}$ if $L(a)=0$ then $a$ is a scalar.
\end{lemma}
\begin{proof}
For any $\gamma\in \mathcal{J}$ by the condition (2) we have
\begin{eqnarray*}
L(\alpha_{\gamma}(a))\le \pa
\dim(\gamma)\overline{\chi_{\gamma}}\pa_1\cdot L(a)=0.
\end{eqnarray*}
By conditions (3) and (4)
we see that $\alpha_{\gamma}(a)=0$ for $\gamma\neq \gamma_0$ and
that $\alpha_{\gamma_0}(a)\in \Ce$. Hence
$\alpha_{\gamma}(a-\alpha_{\gamma_0}(a))=0$ for all $\gamma \in
\hat{G}$. Then Lemma~\ref{comp vanish to vanish:lemma}
tells us that $a=\alpha_{\gamma_0}(a)\in \Ce$.
\end{proof}

\begin{lemma} \label{condition 3*:lemma}
For any $R\ge 0$ the ball
\begin{eqnarray*}
\mathcal{D}_{R}(\mathcal{A})=
\{a\in \mathcal{A}:L(a)\le 1, \pa a\pa \le R\}
\end{eqnarray*}
is totally bounded.
\end{lemma}
\begin{proof}
For any $\epsilon>0$ by Lemma~\ref{finite approx*:lemma}
there is some finite subset $\mathcal{J}\subseteq \hat{G}$ such that
for every $v\in \mathcal{D}_{R}(\mathcal{A})$ there exists $v'\in \mathcal{D}_{R}(\mathcal{A}_{\mathcal{J}})$ with $\pa v-v'\pa<\epsilon$.
Let $M=\max{\{\pa
  \dim(\gamma)\overline{\chi_{\gamma}}\pa_1:\gamma\in \mathcal{J}\}}$.
For any
$a=\sum_{\gamma\in \mathcal{J}}a_{\gamma}\in \mathcal{D}_{R}(\mathcal{A}_{\mathcal{J}})$ and $\gamma\in \mathcal{J}$ we have
\begin{eqnarray*}
\pa a_{\gamma}\pa =\pa \alpha_{\dim(\gamma)\overline{\chi_{\gamma}}}(a)\pa \le
\pa \dim(\gamma)\overline{\chi_{\gamma}}\pa_1\cdot \pa a\pa \le M\cdot R,
\end{eqnarray*}
and by the condition (2)
\begin{eqnarray*}
L(a_{\gamma})=L(\alpha_{\dim(\gamma)\overline{\chi_{\gamma}}}(a))\le
\pa \dim(\gamma)\overline{\chi_{\gamma}}\pa_1 \cdot L(a)\le M.
\end{eqnarray*}
Therefore
\begin{eqnarray*}
\mathcal{D}_{R}(\mathcal{A}_{\mathcal{J}})\subseteq
\{\sum_{\gamma\in \mathcal{J}}a_{\gamma}\in
\mathcal{A}_{\mathcal{J}}:a_{\gamma}\in \mathcal{A}_{\gamma}, \, \,
L(a_{\gamma})\le M,
\, \, \pa a_{\gamma}\pa \le M\cdot R\}.
\end{eqnarray*}
By the conditions (3), (4) and Proposition~\ref{criterion of
Lip:prop}
 the latter set is totally bounded. Then
$\mathcal{D}_{R}(\mathcal{A}_{\mathcal{J}})$ is totally bounded.
Since $\epsilon$ is arbitrary, $\mathcal{D}_{R}(\mathcal{A})$ is also totally bounded.
\end{proof}

\begin{lemma} \label{condition 2*:lemma}
We have
\begin{eqnarray*}
\pa \cdot \pa^{\sim}\le
\Bigl(r_{\mathcal{B}}+C\int_{G}\mathnormal{l}(x)\, dx\Bigr)
L^{\sim}
\end{eqnarray*}
on $\mathcal{A}_{sa}/\Re e$.
\end{lemma}
\begin{proof}
Let $a\in \mathcal{A}_{sa}$ with $L(a)=1$.
Let
$\varphi$ be the constant function $\chi_{\gamma_0}=1$ on $G$.
Then $\alpha_{\varphi}=\alpha_{\gamma_0}$ and $\pa \varphi\pa_1=1$.
As in the proof of Lemma~\ref{finite approx*:lemma} we have $\alpha_{\varphi}(a)\in
(\mathcal{A}^{\alpha})_{sa}$ and
\begin{eqnarray*}
 \pa a-\alpha_{\varphi}(a)\pa
\le
L^{\mathnormal{l}}(a)\int_G\varphi(x)\mathnormal{l}(x)\,
dx\le C\cdot L(a)\int_G\mathnormal{l}(x)\, dx=C\int_G\mathnormal{l}(x)\, dx,
\end{eqnarray*}
where the second inequality comes from the condition (1).
Let $b=\alpha_{\varphi}(a)$.
By the condition (2) we have
\begin{eqnarray*}
L(b)\le \pa \varphi\pa_1\cdot L(a)= 1.
\end{eqnarray*}
Then by Proposition~\ref{criterion of Lip:prop}
\begin{eqnarray*}
r_{\mathcal{B}}\ge \pa \tilde{b}\pa^{\sim}\ge
\pa \tilde{a}\pa^{\sim}-\pa \tilde{a}-\tilde{b}\pa^{\sim}
\ge \pa \tilde{a}\pa^{\sim}-\pa a- \alpha_{\varphi}(a)\pa
\ge \pa \tilde{a}\pa^{\sim}-C\int_G\mathnormal{l}(x)\, dx.
\end{eqnarray*}
Therefore we have $\pa \cdot \pa^{\sim}\le
(r_{\mathcal{B}}+C\int_G\mathnormal{l}(x)\, dx)
L^{\sim}$.
\end{proof}

Now Theorem~\ref{use of action:thm} follows from
Lemmas~\ref{condition 1*:lemma}-\ref{condition 2*:lemma} and
Proposition~\ref{criterion of Lip:prop} immediately.
\end{proof}

\section{Differential Operators and Seminorms}
\label{DiffSN:sec}

In this section we make preparation for our proof of
Theorem~\ref{theta-deform metric:thm}.
In Section~\ref{proof:sec}
we shall verify the conditions in Theorem~\ref{use of action:thm}
for $(C(M_{\theta}), L_{\theta}, \mathbb{T}^n, I\otimes \tau)$.
The seminorm $L^{\mathnormal{l}}_{\theta}$ on
$C(M_{\theta})$ associated to $I\otimes \tau$
is defined in Definition~\ref{L^l L^l_{theta}:def}.
The main difficulty is to verify the condition (1).
We shall see that it is much more convenient to work on the whole
Hilbert space $\mathcal{H}\bar{\otimes}L^2(\mathcal{A}_{\theta})$
instead of $\mathcal{H}_{\theta}$.
So we have to study the corresponding seminorms
$L^D$ and $L^{\mathnormal{l}}$ on $C(M)\otimes \mathcal{A}_{\theta}$
(see Definitions~\ref{L^D:def} and \ref{L^l L^l_{theta}:def}).
We prove the comparison formula for $L^D$ and $L^{\mathnormal{l}}$
first, in (\ref{compare:eq}). Then we relate them to
$L_{\theta}$ and $L^{\mathnormal{l}}_{\theta}$ by
proving (\ref{same norm theta:eq}).
The information about these various seminorms
 is all hidden in differential operators, which
involve mainly the theory of {\LCTVS}.
Subsections~\ref{Diff:sub} and
\ref{SN:sub} are devoted to analyzing these operators.

\subsection{Differential Operators}
\label{Diff:sub}

In this subsection we assume that $M$ is an oriented Riemannian manifold
with an isometric smooth action $\sigma_M$ of $\mathbb{T}^n$. Our aim
is to derive the formulas (\ref{D times I:eq}), (\ref{diff tr:eq}) and
(\ref{diff in tensor:eq}) below.

Let $Cl^{\Ce}M$ be the complexified Clifford algebra
bundle on $M$. Then its space of smooth sections,
$C^{\infty}(M,Cl^{\Ce}M)$,
 is a {\LCA}
containing $C^{\infty}(M)$ as a central subalgebra, and containing
$C^{\infty}(M,TM^{\Ce})$ as a subspace, where $TM^{\Ce}$ is the
complexified tangent bundle. Using the Riemannian metric, we can
identify $TM$ and $T^*M$ canonically. Then $C^{\infty}(M,
T^*M^{\Ce})=C^{\infty}(M, TM^{\Ce})$ is also a subspace of
$C^{\infty}(M,Cl^{\Ce}M)$. Notice that $C^{\infty}(M, S)$ is a
locally convex left module over $C^{\infty}(M,Cl^{\Ce}M)$. Since
$\mathcal{A}^{\infty}_{\theta}$ is nuclear,
the complete tensor products
$C^{\infty}(M)\hat{\otimes}\mathcal{A}^{\infty}_{\theta}$,
$C^{\infty}(M,
TM^{\Ce})\hat{\otimes}\mathcal{A}^{\infty}_{\theta}$ and
$C^{\infty}(M,
T^*M^{\Ce})\hat{\otimes}\mathcal{A}^{\infty}_{\theta}$ can be
thought of as complete injective tensor products, and hence are
are all subspaces of $C^{\infty}(M,Cl^{\Ce}M)
\hat{\otimes}\mathcal{A}^{\infty}_{\theta}$ (see the discussion
after Proposition~\ref{tensor prod of alg:lemma}).

In the same way we think of $C(M, T^*M^{\Ce})=C(M, TM^{\Ce})$ as a
subspace of $C(M,Cl^{\Ce}M)$. Since the $C^*$-algebraic norm on
$Cl^{\Ce}(TM_p)$ extends the inner-product norm on the tangent
space $TM_p$ for each $p\in M$ (see the discussion after
Lemma~\ref{tr:lemma}), clearly the supremum (possibly
$+\infty$-valued) norm on $C(M,Cl^{\Ce}M)$ extends that on $C(M,
TM)$, which is pointwise the inner-product norm.

Clearly the action
of $\mathbb{T}^n$ on the bundle $TM$
extends to an action on the bundle $Cl^{\Ce}M$.
We denote the
induced continuous action on $C^{\infty}(M,Cl^{\Ce}M)$ also by
$\sigma$. Much as in Section~\ref{ConsConnesDubois:sec}, we can
define
\begin{eqnarray*}
C^{\infty}(M_{\theta}, Cl^{\Ce}M)&:=&
(C^{\infty}(M, Cl^{\Ce}M)
\hat{\otimes}\mathcal{A}^{\infty}_{\theta})^{\sigma\hat{\otimes} \tau^{-1}}, \\
C^{\infty}(M_{\theta}, TM^{\Ce})&:=&
(C^{\infty}(M, TM^{\Ce})
\hat{\otimes}\mathcal{A}^{\infty}_{\theta})^{\sigma\hat{\otimes} \tau^{-1}}, \\
C^{\infty}(M_{\theta}, T^*M^{\Ce})&:=& (C^{\infty}(M, T^*M^{\Ce})
\hat{\otimes}\mathcal{A}^{\infty}_{\theta})^{\sigma\hat{\otimes}
  \tau^{-1}}.
\end{eqnarray*}

The differential operator $d:C^{\infty}(M)\rightarrow
C^{\infty}(M, T^*M^{\Ce})$ is a first-order linear operator, and
hence easily seen to be continuous. Then we have the tensor
product linear map $d\hat{\otimes} I:C^{\infty}(M)\hat{\otimes}
\mathcal{A}^{\infty}_{\theta}\rightarrow C^{\infty}(M, T^*M^{\Ce})
\hat{\otimes} \mathcal{A}^{\infty}_{\theta}$. Notice that $d$
commutes with the action $\sigma$. So $d\hat{\otimes} I$ commutes
with $\sigma\hat{\otimes} \tau^{-1}$, and hence maps
$C^{\infty}(M_{\theta})$ into $C^{\infty}(M_{\theta},
T^*M^{\Ce})$. The deformed differential $d_{\theta}$ is then
defined to be  the restriction of $d\hat{\otimes} I$ to
$C^{\infty}(M_{\theta})$.

For any $f\in C^{\infty}(M)$ we have
\begin{eqnarray} \label{D df:eq}
[D, f]=df  \mbox{ as linear maps on } C^{\infty}(M, S),
\end{eqnarray}
where $df\in C^{\infty}(M, T^*M^{\Ce})\subseteq C^{\infty}(M,
Cl^{\Ce}M)$ acts on $C^{\infty}(M, S)$ via the left
$C^{\infty}(M,Cl^{\Ce}M)$-module structure of $C^{\infty}(M, S)$.
Then it is easy to see that for any $f\in C^{\infty}(M)
\otimes_{alg} \mathcal{A}^{\infty}_{\theta}$ we have
\begin{eqnarray*}
[D\otimes I, f]=(d\otimes I)f \mbox{ as linear maps on }
C^{\infty}(M, S)\otimes_{alg} \mathcal{A}^{\infty}_{\theta}.
\end{eqnarray*}
This means that the bilinear maps $(f, \psi)\mapsto
[D\hat{\otimes} I, f](\psi)$ and $(f, \psi)\mapsto
((d\hat{\otimes} I)f)(\psi)$ from $W:=(C^{\infty}(M) \hat{\otimes}
\mathcal{A}^{\infty}_{\theta})\times (C^{\infty}(M,
S)\hat{\otimes} \mathcal{A}^{\infty}_{\theta})$ to $C^{\infty}(M,
S)\hat{\otimes} \mathcal{A}^{\infty}_{\theta}$ coincide on the
dense subspace $(C^{\infty}(M) \otimes_{alg}
\mathcal{A}^{\infty}_{\theta})\times (C^{\infty}(M,
S)\otimes_{alg} \mathcal{A}^{\infty}_{\theta})$. Since both of
them are (jointly) continuous, they coincide on the whole of $W$.
In other words, for any $f\in C^{\infty}(M) \hat{\otimes}
\mathcal{A}^{\infty}_{\theta}$ we have
\begin{eqnarray} \label{D times I:eq}
[D\hat{\otimes} I, f]=(d\hat{\otimes} I)f \mbox{ as linear maps on
} C^{\infty}(M, S)\hat{\otimes} \mathcal{A}^{\infty}_{\theta}.
\end{eqnarray}

The canonical $\Re$-bilinear pairing $C^{\infty}(M, TM)\times
C^{\infty}(M, T^*M)\rightarrow C^{\infty}(M)$ extends to a
$\Ce$-bilinear pairing $C^{\infty}(M, TM^{\Ce})\times
C^{\infty}(M, T^*M^{\Ce})\rightarrow C^{\infty}(M)$, which is
clearly continuous. For any $Y\in C^{\infty}(M, TM^{\Ce})$ let
$\mathfrak{i}_Y$ be the corresponding contraction $C^{\infty}(M,
T^*M^{\Ce})\rightarrow C^{\infty}(M)$. Then we have the
tensor-product map $\mathfrak{i}_Y\hat{\otimes}I:C^{\infty}(M,
T^*M^{\Ce})\hat{\otimes} \mathcal{A}^{\infty}_{\theta} \rightarrow
C^{\infty}(M)\hat{\otimes}$ $ \mathcal{A}^{\infty}_{\theta}$. Let
$\partial_Y:C^{\infty}(M)\rightarrow C^{\infty}(M)$ be the
derivation with respect to $Y$. Since $\partial_Y$ is a
first-order linear operator, it is continuous. Then we also have
the tensor-product map
$\partial_Y\hat{\otimes}I:C^{\infty}(M)\hat{\otimes}
\mathcal{A}^{\infty}_{\theta} \rightarrow
C^{\infty}(M)\hat{\otimes} \mathcal{A}^{\infty}_{\theta}$. For any
$f\in C^{\infty}(M)$ it is trivial to see that
\begin{eqnarray*}
\partial_Y(f)=\mathfrak{i}_Y(df).
\end{eqnarray*}
Then for any $f\in C^{\infty}(M)\otimes_{alg}
\mathcal{A}^{\infty}_{\theta}$ clearly
\begin{eqnarray*}
(\partial_Y\otimes I)(f)=((\mathfrak{i}_Y\otimes I)\circ (d\otimes I))(f).
\end{eqnarray*}
By the same argument as for (\ref{D times I:eq}), for any $f\in
C^{\infty}(M)\hat{\otimes}\mathcal{A}^{\infty}_{\theta}$ we then
have
\begin{eqnarray} \label{diff:eq}
(\partial_Y\hat{\otimes} I)(f)=
((\mathfrak{i}_Y\hat{\otimes} I)\circ (d\hat{\otimes} I))(f).
\end{eqnarray}

Since the tracial state $tr:Cl^{\Ce}(TM_p)\rightarrow \Ce$ in
Lemma~\ref{tr:lemma} is invariant under the action of $SO(TM_p)$
for each $p\in M$, we can use them pointwisely to define a linear
map $C^{\infty}(M, Cl^{\Ce}M)\to C^{\infty}(M)$, which is clearly
continuous. We denote this map also by $tr$. Then $tr$ is still
tracial in the sense that $tr(f\cdot g)=tr(g\cdot f)$ for any
$f,g\in C^{\infty}(M, Cl^{\Ce}M)$. We have the tensor-product
linear map $tr\hat{\otimes}I:C^{\infty}(M, Cl^{\Ce}M)\hat{\otimes}
\mathcal{A}^{\infty}_{\theta}\rightarrow
C^{\infty}(M)\hat{\otimes}\mathcal{A}^{\infty}_{\theta}$. For any
$Y\in C^{\infty}(M, TM^{\Ce})\subseteq C^{\infty}(M, Cl^{\Ce}M)$
and $Z\in C^{\infty}(M, T^*M^{\Ce})\subseteq  C^{\infty}(M,
Cl^{\Ce}M)$, recalling that we have a canonical identification of
$ C^{\infty}(M, TM^{\Ce})$ and $C^{\infty}(M, T^*M^{\Ce})$, we get
\begin{eqnarray*}
tr(Y\cdot Z)=\frac{1}{2}tr(Y\cdot Z+Z\cdot
Y)=\frac{1}{2}tr(-2<Y,Z>)=-<Y, Z>=-\mathfrak{i}_Y(Z),
\end{eqnarray*}
where $Y\cdot Z$ is the multiplication in $C^{\infty}(M,
Cl^{\Ce}M)$, and $<\cdot , \cdot>$ is the $C^{\infty}(M)$-valued
$C^{\infty}(M)$-bilinear pairing on $C^{\infty}(M, TM^{\Ce})$. So
$\mathfrak{i}_Y=tr\circ (-Y)$ on \\$C^{\infty}(M, T^*M^{\Ce})$.
Then $\mathfrak{i}_Y\otimes I=(tr\otimes I)\circ ((-Y)\otimes 1)$
on $C^{\infty}(M,
T^*M^{\Ce})\otimes_{alg}\mathcal{A}^{\infty}_{\theta}$. Since both
$\mathfrak{i}_Y\hat{\otimes} I$ and $(tr\hat{\otimes}I)\circ
((-Y)\otimes 1)$ are continuous maps from \\$C^{\infty}(M,
T^*M^{\Ce})\hat{\otimes} \mathcal{A}^{\infty}_{\theta}$ to
$C^{\infty}(M)\hat{\otimes} \mathcal{A}^{\infty}_{\theta}$, we get
\begin{eqnarray} \label{tr:eq}
\mathfrak{i}_Y\hat{\otimes} I=(tr\hat{\otimes}I)\circ
((-Y)\otimes 1).
\end{eqnarray}
as maps $C^{\infty}(M, T^*M^{\Ce})\hat{\otimes}
\mathcal{A}^{\infty}_{\theta}\rightarrow
C^{\infty}(M)\hat{\otimes} \mathcal{A}^{\infty}_{\theta}$.
Combining (\ref{diff:eq}) and (\ref{tr:eq}) together, for any
$f\in C^{\infty}(M)\hat{\otimes}\mathcal{A}^{\infty}_{\theta}$ we
get
\begin{eqnarray} \label{diff tr:eq}
(\partial_Y\hat{\otimes} I)(f)=
((tr\hat{\otimes}I)\circ
((-Y)\otimes 1)\circ (d\hat{\otimes} I))(f).
\end{eqnarray}

Let  $\Lie(\mathbb{T}^n)$ be the Lie algebra of $\mathbb{T}^n$.
For any $X\in \Lie(\mathbb{T}^n)$ we denote by $X^{\#}$
the vector field on $M$ generated by $X$.

\begin{lemma} \label{diff in tensor:lemma}
For any $X\in \Lie(\mathbb{T}^n)$ and any $ f\in
C^{\infty}(M)\hat{\otimes}\mathcal{A}^{\infty}_{\theta}$ we have
\begin{eqnarray} \label{diff in tensor:eq}
\lim_{t\to 0}\frac{(\sigma_{e^{tX}}\hat{\otimes} I)(f)-f}{t}=
(\partial_{-X^{\#}}\hat{\otimes} I)(f).
\end{eqnarray}
\end{lemma}
\begin{proof}
For any $f\in C^{\infty}(M)$ and $x\in \mathbb{T}^n$
 clearly
\begin{eqnarray*}
(\partial_{-X^{\#}})(f)&=&\lim_{t\to 0}
\frac{\sigma_{e^{tX}}(f)-f}{t},\\
(\partial_{-X^{\#}})(\sigma_x(f))&=&\lim_{t\to 0}
\frac{\sigma_{e^{tX}}(\sigma_x(f))-\sigma_x(f)}{t}\\
&=& \lim_{t\to 0}
\sigma_x(\frac{\sigma_{e^{tX}}(f)-f}{t})=
\sigma_x(\partial_{-X^{\#}}(f)),
\end{eqnarray*}
where the limits are taken with respect to the locally convex topology
in $ C^{\infty}(M)$.
(Here we have $-X^{\#}$ instead of $X^{\#}$ in the
first equation because $(\sigma_{e^{tX}}(f))(p)=f(\sigma_{e^{-tX}}(p))$
for any $p\in M$.)
So we see that the map
$t\mapsto \partial_{-X^{\#}}(\sigma_{e^{tX}}(f))$ is continuous.
When $M$ is compact, we know that
\begin{eqnarray} \label{int:eq}
\sigma_{e^{tX}}(f)-f=\int^t_0
\partial_{-X^{\#}}(\sigma_{e^{sX}}(f))\, ds=\int^t_0
\sigma_{e^{sX}}(\partial_{-X^{\#}}(f))\, ds,
\end{eqnarray}
where the integral is taken with respect to the supremum norm topology in
$C(M)$. Notice that the inclusion $C^{\infty}(M)\hookrightarrow C(M)$ is
continuous when $C^{\infty}(M)$ is endowed with the locally convex
topology and $C(M)$ is endowed with the norm topology. By
Proposition~\ref{linear of integ:prop} the integral $\int^t_0
\sigma_{e^{sX}}(\partial_{-X^{\#}}(f))\, ds$ is also defined in
$C^{\infty}(M)$, and is mapped to the corresponding integral in $C(M)$
under the inclusion $C^{\infty}(M)\hookrightarrow C(M)$. Therefore we see
that (\ref{int:eq}) also holds with respect to the locally convex topology
in $C^{\infty}(M)$. For noncompact $M$, since the locally convex
topology on $C^{\infty}(M)$ is defined using seminorms from compact subsets of
local trivializations, it is
easy to see that (\ref{int:eq}) still holds.

Now for any $f\in
C^{\infty}(M)\otimes_{alg}\mathcal{A}^{\infty}_{\theta}$ clearly
we have
\begin{eqnarray*}
(\sigma_{e^{tX}}\otimes I)(f)-f=\int^t_0
(\sigma_{e^{sX}}\otimes I)((\partial_{-X^{\#}}\otimes I)(f))\, ds
\end{eqnarray*}
in $ C^{\infty}(M)\hat{\otimes}\mathcal{A}^{\infty}_{\theta}$. For
fixed $X$ notice that $f\mapsto (\sigma_{e^{tX}}\hat{\otimes}
I)(f)-f$ is a continuous map from
$C^{\infty}(M)\hat{\otimes}\mathcal{A}^{\infty}_{\theta}$ to
itself. It is also easy to see that both $f\mapsto
(\partial_{-X^{\#}}\hat{\otimes} I)(f)$ and $f\mapsto \int^t_0
(\sigma_{e^{sX}}\hat{\otimes} I)(f)\, ds$ are continuous maps from
$C^{\infty}(M)\hat{\otimes}\mathcal{A}^{\infty}_{\theta}$ to
itself. So the map $f\mapsto \int^t_0
(\sigma_{e^{sX}}\hat{\otimes} I) ((\partial_{-X^{\#}}\hat{\otimes}
I)(f))\, ds$
 from $C^{\infty}(M)\hat{\otimes}\mathcal{A}^{\infty}_{\theta}$ to
itself is continuous. Therefore, for any $f\in
C^{\infty}(M)\hat{\otimes}\mathcal{A}^{\infty}_{\theta}$ we have
\begin{eqnarray*}
(\sigma_{e^{tX}}\hat{\otimes} I)(f)-f=\int^t_0
(\sigma_{e^{sX}}\hat{\otimes} I)
((\partial_{-X^{\#}}\hat{\otimes} I)(f))\, ds.
\end{eqnarray*}
Now (\ref{diff in tensor:eq}) follows from
Proposition~\ref{integ to diff:prop}.
\end{proof}

\subsection{Seminorms}
\label{SN:sub}

In this subsection we assume that $M$ is an $m$-dimensional
compact Spin manifold, and that the action $\sigma_M$ lifts to an
action on $S$. Notice that the fibres of  $Cl^{\Ce}M$ are all
isomorphic to the $C^*$-algebra $Cl^{\Ce}(\Re^m)$, where $\Re^m$
is the standard $m$-dimensional Euclidean space. Clearly
$C^{\infty}(M, Cl^{\Ce}M)$ generates a continuous field of
$C^*$-algebras \cite[Secton 10.3]{Dixmier77} over $M$ with
continuous sections $\Gamma'=C(M, Cl^{\Ce}M)$. Recall that
$\mathcal{H}$ is the Hilbert space completion of $C^{\infty}(M,
S)$. So the algebra $C(M,Cl^{\Ce}M)$ has a natural faithful
representation on $\mathcal{H}$. It is easy to see that the
inclusion $C^{\infty}(M,Cl^{\Ce}M)\hookrightarrow C(M,Cl^{\Ce}M)$
is continuous with respect to the locally convex topology on
$C^{\infty}(M,Cl^{\Ce}M)$ and the norm topology on
$C(M,Cl^{\Ce}M)$. Just as in the case of
$C^{\infty}(M)\hat{\otimes}\mathcal{A}^{\infty}_{\theta}
\rightarrow C(M)\otimes \mathcal{A}_{\theta}$ in
Section~\ref{ConsConnesDubois:sec}, we have a
$\mathbb{T}^n$-equivariant continuous linear map $C^{\infty}(M,
Cl^{\Ce}M) \hat{\otimes}\mathcal{A}^{\infty}_{\theta} \rightarrow
C(M,Cl^{\Ce}M)\otimes \mathcal{A}_{\theta}$ extending this former
one. We still denote it by $\Psi$. As in Lemmas~\ref{injective
Phi:lemma} and \ref{injective Psi:lemma}, $\Psi$ is in fact
injective. Clearly $\Psi$ is a $*$-algebra homomorphism. Let
$C(M_{\theta}, Cl^{\Ce}M_{\theta})$ be $(C(M, Cl^{\Ce}M)\otimes
\mathcal{A}_{\theta})^{\sigma\otimes \tau^{-1}}$. We also have the
homomorphism $C^{\infty}(M_{\theta},Cl^{\Ce}M_{\theta})\rightarrow
B(\mathcal{H}_{\theta})$, which we still denote by
$\Psi_{\theta}$.

\begin{proposition} \label{commutator psi:prop}
For any $f\in C^{\infty}(M)\hat{\otimes}
\mathcal{A}^{\infty}_{\theta}$
 the domain of $D^{L^2}$ is stable under $\Psi(f)$, and
\begin{eqnarray} \label{commutator psi:eq}
[D^{L^2}, \Psi(f)] =
\Psi((d\hat{\otimes} I)f).
\end{eqnarray}
When $f$ is in $C^{\infty}(M_{\theta})$,
 the domain of $D^{L^2}_{\theta}$ is stable under $\Psi_{\theta}(f)$, and
\begin{eqnarray} \label{commutator psi theta:eq}
[D^{L^2}_{\theta}, \Psi_{\theta}(f)] =
\Psi_{\theta}(d_{\theta}f).
\end{eqnarray}
\end{proposition}
\begin{proof}
By Lemma~\ref{tensor prod of alg:lemma} $C^{\infty}(M,
S)\hat{\otimes} \mathcal{A}^{\infty}_{\theta}$ is a locally convex
left module over the algebra $C^{\infty}(M,Cl^{\Ce}M)$
$\hat{\otimes} \mathcal{A}^{\infty}_{\theta}$. So we have the
continuous maps:
\begin{eqnarray*}
(C^{\infty}(M,Cl^{\Ce}M)\hat{\otimes}
\mathcal{A}^{\infty}_{\theta})\times (C^{\infty}(M,
S)\hat{\otimes} \mathcal{A}^{\infty}_{\theta})\rightarrow
C^{\infty}(M, S)\hat{\otimes}
\mathcal{A}^{\infty}_{\theta}\overset{\Phi}\rightarrow
\mathcal{H}\bar{\otimes} L^2(\mathcal{A}_{\theta}).
\end{eqnarray*}
On the other hand, we have continuous maps:
\begin{eqnarray*}
& & (C^{\infty}(M,Cl^{\Ce}M)\hat{\otimes}
\mathcal{A}^{\infty}_{\theta})\times (C^{\infty}(M,
S)\hat{\otimes} \mathcal{A}^{\infty}_{\theta})
\\& & \overset{\Psi\times \Phi}\longrightarrow
B(\mathcal{H}\bar{\otimes}
  L^2(\mathcal{A}_{\theta}))
\times\mathcal{H}\bar{\otimes} L^2(\mathcal{A}_{\theta}) \rightarrow
\mathcal{H}\bar{\otimes} L^2(\mathcal{A}_{\theta}).
\end{eqnarray*}
The two compositions coincide on
$(C^{\infty}(M,Cl^{\Ce}M)\otimes_{alg}
\mathcal{A}^{\infty}_{\theta})\times (C^{\infty}(M,
S)\otimes_{alg} \mathcal{A}^{\infty}_{\theta})$. So they coincide
on the whole of $(C^{\infty}(M,Cl^{\Ce}M)\hat{\otimes}
\mathcal{A}^{\infty}_{\theta})\times (C^{\infty}(M,
S)\hat{\otimes}$ $\mathcal{A}^{\infty}_{\theta})$. In other words,
for any $f\in C^{\infty}(M,Cl^{\Ce}M)\hat{\otimes}
\mathcal{A}^{\infty}_{\theta}$ and any $\psi \in C^{\infty}(M,
S)\hat{\otimes}$ $\mathcal{A}^{\infty}_{\theta}$ we have
\begin{eqnarray} \label{module hom:eq}
\Psi(f)\cdot \Phi(\psi)=\Phi(f\psi).
\end{eqnarray}
Then for any $f\in C^{\infty}(M,Cl^{\Ce}M)\hat{\otimes}
\mathcal{A}^{\infty}_{\theta}$ and $\psi \in C^{\infty}(M,
S)\hat{\otimes} \mathcal{A}^{\infty}_{\theta}$ we have
\begin{eqnarray*}
& &\Phi([D\hat{\otimes}I, f](\psi))\\&=&
\Phi((D\hat{\otimes}I)(f\psi)-f((D\hat{\otimes}I)\psi))
\overset{(\ref{D Psi:eq})}= D^{L^2}(\Phi(f\psi))-
\Phi(f((D\hat{\otimes}I)\psi))  \\
&\overset{(\ref{module hom:eq})}=& D^{L^2}((\Psi(f))(\Phi(\psi)))-
\Psi(f)\cdot \Phi((D\hat{\otimes}I)\psi)
\\&\overset{(\ref{D Psi:eq})}=&
D^{L^2}((\Psi(f))(\Phi(\psi)))-\Psi(f)(D^{L^2}(\Phi(\psi)))
= [D^{L^2}, \Psi(f)](\Phi(\psi)).
\end{eqnarray*}
So for any $f\in C^{\infty}(M,Cl^{\Ce}M)\hat{\otimes}
\mathcal{A}^{\infty}_{\theta}$ we have
\begin{eqnarray} \label{exchange D:eq}
\Phi\circ [D\hat{\otimes}I, f]= [D^{L^2}, \Psi(f)]\circ \Phi
\end{eqnarray}
as linear maps from $C^{\infty}(M,
S)\hat{\otimes}\mathcal{A}^{\infty}_{\theta}$ to
$\mathcal{H}\bar{\otimes} L^2(\mathcal{A}_{\theta})$. When $f$ is
in $C^{\infty}(M)\hat{\otimes} \mathcal{A}^{\infty}_{\theta}$ we
also have
\begin{eqnarray*}
\Phi\circ [D\hat{\otimes}I, f]
\overset{(\ref{D times I:eq})}=\Phi\circ ((d\hat{\otimes} I)f)
\overset{(\ref{module hom:eq})}=\Psi((d\hat{\otimes} I)f)\circ \Phi.
\end{eqnarray*}
Therefore, for any $f\in C^{\infty}(M)\hat{\otimes}
\mathcal{A}^{\infty}_{\theta}$ we have
\begin{eqnarray} \label{D Psi Phi:eq}
[D^{L^2}, \Psi(f)]\circ \Phi =
\Psi((d\hat{\otimes} I)f)\circ \Phi.
\end{eqnarray}
For any $z$ in the domain of $D^{L^2}$ take a net $\psi_j$ in
$C^{\infty}(M, S)\hat{\otimes}\mathcal{A}^{\infty}_{\theta}$ with
$\Phi(\psi_j)\to z$ and $D^{L^2}(\Phi(\psi_j))\to D^{L^2}(z)$.
Then
\begin{eqnarray*}
D^{L^2}((\Psi(f))(\Phi(\psi_j)))&\overset{(\ref{D Psi Phi:eq})}=&
(\Psi(f))(D^{L^2}(\Phi(\psi_j)))+
\Psi((d\hat{\otimes}I)(f))(\Phi(\psi_j))\\
&\to &(\Psi(f))(D^{L^2}(z))+
\Psi((d\hat{\otimes}I)(f))(z),
\end{eqnarray*}
and
\begin{eqnarray*}
(\Psi(f))(\Phi(\psi_j))\to (\Psi(f))(z).
\end{eqnarray*}
So $(\Psi(f))(z)$ is in the domain of $D^{L^2}$, and
\begin{eqnarray*}
D^{L^2}((\Psi(f))(z))=(\Psi(f))(D^{L^2}(z))+
\Psi((d\hat{\otimes}I)(f))(z).
\end{eqnarray*}
Therefore the domain of $D^{L^2}$ is stable under $\Psi(f)$,
 and $[D^{L^2}, \Psi(f)] =
\Psi((d\hat{\otimes} I)f)$.

The assertions about $C^{\infty}(M_{\theta})$ follow from those
about $C^{\infty}(M)\hat{\otimes} \mathcal{A}^{\infty}_{\theta}$.
\end{proof}

By Proposition~\ref{commutator psi:prop} we see that the
commutator $[D^{L^2}, f]$ is bounded for any $f\in
\Psi(C^{\infty}(M)\hat{\otimes} \mathcal{A}^{\infty}_{\theta})$.
Corresponding to $L_{\theta}$ defined in
Definition~\ref{L_{theta}:def} we have:
\begin{definition} \label{L^D:def}
We define a seminorm, denoted by $L^D$, on $C(M)\otimes
\mathcal{A}_{\theta}$ by
\begin{eqnarray*}
L^D(f):=\begin{cases}
\pa [D^{L^2}, f]\pa, & \mbox{ if } f\in
\Psi(C^{\infty}(M)\hat{\otimes}
\mathcal{A}^{\infty}_{\theta});\\
+\infty, & \mbox{ otherwise }.
\end{cases}
\end{eqnarray*}
\end{definition}

Fix an inner product on $\Lie(\mathbb{T}^n)$, and use it to get a
translation-invariant Riemannian metric on $\mathbb{T}^n$ in the usual way.
We get
a length function
$\mathnormal{l}$ on $\mathbb{T}^n$ by
setting $\mathnormal{l}(x)$ to be the geodesic distance from $x$ to $e_{\mathbb{T}^n}$
for $x\in \mathbb{T}^n$.
Notice that $I\otimes \tau=\sigma\otimes I$ is a nontrivial action of $\mathbb{T}^n$ on
$C(M_{\theta})$. To make use of Theorem~\ref{use of action:thm} we
define two seminorms:
\begin{definition} \label{L^l L^l_{theta}:def}
We define a (possibly $+\infty$-valued) seminorm $L^{\mathnormal{l}}$ on
$C(M)\otimes \mathcal{A}_{\theta}$ for the action $\sigma\otimes I$ via
 (\ref{def of L:eq}):
\begin{eqnarray*}
L^{\mathnormal{l}}(f):=
\sup\{ \frac{\pa (\sigma\otimes I)_x(f)-f\pa}{\mathnormal{l}(x)}:
 x\in  \mathbb{T}^n,  x\neq e_{
  \mathbb{T}^n} \}.
\end{eqnarray*}
We also define a (possibly $+\infty$-valued) seminorm $L^{\mathnormal{l}}_{\theta}$
on $C(M_{\theta})$
 for the action $I\otimes \tau$:
\begin{eqnarray*}
L^{\mathnormal{l}}_{\theta}(f):=
\sup\{ \frac{\pa (I\otimes \tau)_x(f)-f\pa}{\mathnormal{l}(x)}:
 x\in  \mathbb{T}^n, x\neq e_{
  \mathbb{T}^n} \}.
\end{eqnarray*}
\end{definition}
Then
\begin{eqnarray} \label{Lltheta to Ll:eq}
L^{\mathnormal{l}}_{\theta}=L^{\mathnormal{l}}
\end{eqnarray}
on $C(M_{\theta})$, because there $I\otimes \tau=\sigma\otimes I$.

Our  first  key technical fact is the following comparison between
$L^{\mathnormal{l}}$ and $L^D$:

\begin{proposition} \label{compare:prop}
Let $C$ be the norm of the linear map $X\mapsto X^{\#}$
from $\Lie(\mathbb{T}^n)$ to $C^{\infty}(M, TM)\subseteq
C(M, Cl^{\Ce}M)$. Then on $C(M)\otimes \mathcal{A}_{\theta}$ we have
\begin{eqnarray} \label{compare:eq}
L^{\mathnormal{l}}\le C\cdot L^D.
\end{eqnarray}
\end{proposition}
\begin{proof}
Let $X\in \Lie(\mathbb{T}^n)$. For any $f\in
C^{\infty}(M)\hat{\otimes} \mathcal{A}^{\infty}_{\theta}$ we have
\begin{eqnarray*}
(\Psi\circ (\partial_{-X^{\#}}\hat{\otimes} I))(f)
&\overset{(\ref{diff in tensor:eq})}=&
\Psi(\lim_{t\to 0}\frac{(\sigma_{e^{tX}}\hat{\otimes}
  I)(f)-f}{t})\\
&=& \lim_{t\to 0}\Psi(\frac{(\sigma_{e^{tX}}\hat{\otimes}
  I)(f)-f}{t})\\ &=& \lim_{t\to 0}\frac{(\sigma_{e^{tX}}\otimes
  I)(\Psi(f))-\Psi(f)}{t}.
\end{eqnarray*}
It follows immediately that $\Psi(f)$ is once-differentiable with
respect to the action $\sigma\otimes I$. In fact, $\Psi(f)$ is
easily seen to be smooth for the action $\sigma\otimes I$, though
we don't need this fact here. By \cite[Proposition 8.6]{Rieffel00}
\begin{eqnarray*}
L^{\mathnormal{l}}(\Psi(f))=\sup_{\pa X\pa=1}\pa
\lim_{t\to 0}\frac{(\sigma_{e^{tX}}\otimes
  I)(\Psi(f))-\Psi(f)}{t} \pa.
\end{eqnarray*}
Then we get
\begin{eqnarray*}
L^{\mathnormal{l}}(\Psi(f))&=&\sup_{\pa X\pa=1}\pa
(\Psi\circ (\partial_{-X^{\#}}\hat{\otimes} I))(f) \pa \\
&\overset{(\ref{diff tr:eq})}=& \sup_{\pa X\pa=1}\pa
(\Psi\circ
(tr\hat{\otimes}I)\circ
((-X^{\#})\otimes 1)\circ (d\hat{\otimes} I))(f) \pa.
\end{eqnarray*}
Notice that
the linear map
$tr:C^{\infty}(M,Cl^{\Ce}M)\rightarrow C^{\infty}(M)$
extends to $C(M,Cl^{\Ce}M)\\ \rightarrow C(M)$, which we still
denote by $tr$.
By Lemma~\ref{tr:lemma} the map
$tr:Cl^{\Ce}(\Re^m)\rightarrow \Ce$ is positive.
Then so is
$tr:C(M,Cl^{\Ce}M)\rightarrow C(M)$.
Since $C(M)$ is commutative, $tr:C(M,Cl^{\Ce}M)\rightarrow C(M)$
 is completely positive \cite[Lemma 5.1.4]{ER00}.
Then we have the tensor-product completely positive map
\cite[Proposition 8.2]{Lance95}
$tr\otimes I:
C(M,Cl^{\Ce}M)\otimes \mathcal{A}_{\theta}\rightarrow
C(M)\otimes \mathcal{A}_{\theta}$.
Consequently, we have $\pa tr\otimes I\pa=\pa (tr\otimes I)(1\otimes
1)\pa=1$ \cite[Lemma 5.1.1]{ER00}.
In fact, $tr\otimes I$ is easily seen to be a conditional
expectation in the sense of \cite[Exercise 8.7.23]{KR97b},
though we don't need this fact here.
Clearly
\begin{eqnarray} \label{tr Psi:eq}
(tr\otimes I)\circ \Psi=\Psi\circ (tr\hat{\otimes}I)
\end{eqnarray}
holds on $C^{\infty}(M, Cl^{\Ce}M)\otimes_{alg}
\mathcal{A}^{\infty}_{\theta}$. Since both maps here are
continuous, (\ref{tr Psi:eq}) holds on the whole of $C^{\infty}(M,
Cl^{\Ce}M)\hat{\otimes} \mathcal{A}^{\infty}_{\theta}$. For any
$Y\in C^{\infty}(M, Cl^{\Ce}M)\subseteq C(M,Cl^{\Ce}M)$, we have
\begin{eqnarray*}
& &\pa (\Psi\circ
(tr\hat{\otimes}I)\circ
((-Y)\otimes 1)\circ (d\hat{\otimes} I))(f) \pa \\
&\overset{(\ref{tr Psi:eq})}=& \pa
((tr\otimes I)\circ \Psi \circ
((-Y)\otimes 1)\circ (d\hat{\otimes} I))(f) \pa \\
&\le & \pa (\Psi \circ
((-Y)\otimes 1)\circ (d\hat{\otimes} I))(f) \pa
= \pa
\Psi( (-Y)\otimes 1)\cdot \Psi( (d\hat{\otimes} I)(f)) \pa \\
&\le & \pa
\Psi( (-Y)\otimes 1)\pa \cdot \pa \Psi( (d\hat{\otimes} I)(f))
\pa
=\pa Y
\pa \cdot \pa \Psi( (d\hat{\otimes} I)(f))
\pa.
\end{eqnarray*}
Recall that $X^{\#}\in
C^{\infty}(M, TM)\subseteq C^{\infty}(M, Cl^{\Ce}M)$.
Therefore
\begin{eqnarray*}
L^{\mathnormal{l}}(\Psi(f))
&=& \sup_{\pa X\pa=1}\pa
(\Psi\circ
(tr\hat{\otimes}I)\circ
((-X^{\#})\otimes 1)\circ (d\hat{\otimes} I))(f) \pa \\
&\le &\sup_{\pa X\pa=1}\pa X^{\#}
\pa \cdot \pa \Psi( (d\hat{\otimes} I)(f))
\pa = C \pa \Psi( (d\hat{\otimes} I)(f))\pa \\
&\overset{(\ref{commutator psi:eq})}=&
C\pa [D^{L^2}, \Psi(f)]\pa
=C\cdot L^D(\Psi(f))
\end{eqnarray*}
as desired.
\end{proof}

\subsection{Restriction Map}
\label{RM:sub}

Our goal in this subsection is to prove the second key technical fact:

\begin{proposition} \label{same norm theta:prop}
The restriction map from $C(M_{\theta},Cl^{\Ce}M)$ to
$B(\mathcal{H}_{\theta})$ is isometric. In particular,
for any $f\in C^{\infty}(M_{\theta},Cl^{\Ce}M)$ we have
\begin{eqnarray} \label{same norm theta:eq}
\pa \Psi(f)\pa=\pa \Psi_{\theta}(f)\pa.
\end{eqnarray}
\end{proposition}

First of all,  Proposition~\ref{same norm theta:prop} justifies our
way of taking $C(M_{\theta})$ as a subalgebra of $B(\mathcal{H}_{\theta})$ via
restriction to $\mathcal{H}_{\theta}$. Secondly,
it enables us to compute $L_{\theta}$ using
our seminorm $L^D$ in Subsection~\ref{SN:sub}, and hence to compare it with
$L^{\mathnormal{l}}_{\theta}$:
\begin{corollary} \label{Ltheta to LD:coro}
On $C(M_{\theta})$ we have
\begin{eqnarray} \label{Ltheta to LD:eq}
L_{\theta}=L^D,
\end{eqnarray}
and
\begin{eqnarray} \label{Lltheta to Ltheta:eq}
L^{\mathnormal{l}}_{\theta}\le C\cdot L_{\theta}.
\end{eqnarray}
\end{corollary}
\begin{proof}
We prove (\ref{Ltheta to LD:eq}) first.
Since $\Psi$ is injective it suffices to show (\ref{Ltheta to LD:eq}) on
$\Psi(C^{\infty}(M_{\theta}))$.
For any $f\in C^{\infty}(M_{\theta})$
we have
\begin{eqnarray*}
L^D(\Psi(f))&=& \pa [D^{L^2}, \Psi(f)]\pa
\overset{(\ref{commutator psi:eq})}= \pa \Psi((d\hat{\otimes} I)f)\pa
= \pa \Psi(d_{\theta}f)\pa
\overset{(\ref{same norm theta:eq})}= \pa \Psi_{\theta}(d_{\theta}f)\pa \\
&\overset{(\ref{commutator psi theta:eq})}=&
\pa [D^{L^2}_{\theta}, \Psi_{\theta}(f)]\pa
= L_{\theta}(\Psi_{\theta}(f)),
\end{eqnarray*}
which yields (\ref{Ltheta to LD:eq}).
Then on $C(M_{\theta})$ we have
\begin{eqnarray*}
L^{\mathnormal{l}}_{\theta}
\overset{(\ref{Lltheta to Ll:eq})}=L^{\mathnormal{l}}
\overset{(\ref{compare:eq})}\le C\cdot L^D
\overset{(\ref{Ltheta to LD:eq})}=C\cdot L_{\theta}.
\end{eqnarray*}
\end{proof}

Instead of proving Proposition~\ref{same norm theta:prop}
directly, we shall prove  a slightly more general form. Let
$\mathcal{A}$ be a unital $C^*$-algebra with a strongly continuous
action $\sigma$ of $\mathbb{T}^n$, which we shall set to be $C(M,
Cl^{\Ce}M)$ later. Assume that $\mathcal{A}\subseteq
B(\mathcal{H})$ and that $\mathbb{T}^n$ has a strongly continuous
unitary representation on $\mathcal{H}$, which we still denote by
$\sigma$, such that the action $\sigma$ on $\mathcal{A}$ is
induced by conjugation. Then the $C^*$-algebraic spatial tensor
product $\mathcal{A}\otimes \mathcal{A}_{\theta}$ \cite[Appendix
T.5]{Wegge93} acts on $\mathcal{H}\bar{\otimes}
L^2(\mathcal{A}_{\theta})$ faithfully. For any $q\in
\Ze^n=\widehat{\mathbb{T}^n}$ let $(\mathcal{H}\bar{\otimes}
L^2(\mathcal{A}_{\theta}))_q$  be the $q$-isotypic subspace of
$\mathcal{H}\bar{\otimes}
  L^2(\mathcal{A}_{\theta})$
for the action $\sigma\bar{\otimes} \tau^{-1}$.
Notice that
$(\mathcal{H}\bar{\otimes} L^2(\mathcal{A}_{\theta}))_q$
is stable under the action of
$(\mathcal{A}\otimes \mathcal{A}_{\theta})^{\sigma\otimes
  \tau^{-1}}$
for each $q\in \Ze^n$.

\begin{proposition} \label{same norm:prop}
For any $f\in (\mathcal{A}\otimes \mathcal{A}_{\theta})^{\sigma\otimes
  \tau^{-1}}$ and $q\in \Ze^n$ we have
\begin{eqnarray} \label{same norm:eq}
\pa f\pa =\pa f|_{(\mathcal{H}\bar{\otimes}
    L^2(\mathcal{A}_{\theta}))_q}\pa,
\end{eqnarray}
where $(\mathcal{H}\bar{\otimes}
    L^2(\mathcal{A}_{\theta}))_q$ is the $q$-isotypic component of
$\mathcal{H}\bar{\otimes}
    L^2(\mathcal{A}_{\theta})$ under $\sigma\bar{\otimes}\tau^{-1}$.
\end{proposition}
\begin{proof} Think of $-\theta
  q$ as an element of  $\mathbb{T}^n$ via the natural projection
$\Re^n\rightarrow \Re^n/\Ze^n=\mathbb{T}^n$.
For any $p\in \Ze^n$, recalling the skew-symmetric bicharacter
$\omega_{\theta}$ in Section~\ref{Prelim:sec},
we have
\begin{eqnarray*}
u_qu_pu_{-q}=
\omega_{\theta}(q,p)\omega_{\theta}(q+p,-q)u_p=
\omega_{\theta}(q,2p)u_p=\left<p,-\theta q\right>u_p=
\tau_{-\theta q}(u_p).
\end{eqnarray*}
It follows immediately that for any $b\in \mathcal{A}_{\theta}$ we
have
\begin{eqnarray*}
u_qbu_{-q}=\tau_{-\theta q}(b).
\end{eqnarray*}
Consequently, for any $f\in \mathcal{A}\otimes \mathcal{A}_{\theta}$ we have
\begin{eqnarray*}
(1\otimes u_q)f(1\otimes u_{-q})=(I\otimes \tau)_{-\theta q}(f).
\end{eqnarray*}
Therefore
\begin{eqnarray} \label{conj:eq}
(1\otimes u_q)\circ f\circ (1\otimes u_{-q})
=(I\bar{\otimes} \tau)_{-\theta
  q}\circ f\circ (I\bar{\otimes} \tau)_{\theta
  q}
\end{eqnarray}
on $\mathcal{H}\bar{\otimes}L^2(\mathcal{A}_{\theta})$.
Clearly $1\otimes u_{-q}$ is in the $q$-isotypic component of
$\mathcal{A}\otimes \mathcal{A}_{\theta}$ under $\sigma\otimes
\tau^{-1}$. So $1\otimes u_{-q}$ restricted to $(\mathcal{H}\bar{\otimes}
    L^2(\mathcal{A}_{\theta}))_0$ is a unitary map onto
$(\mathcal{H}\bar{\otimes} L^2(\mathcal{A}_{\theta}))_q$.
Since $I\bar{\otimes} \tau$ and
$\sigma\bar{\otimes} \tau^{-1}$ commute with each other,
$I\bar{\otimes}\tau$ preserves
$(\mathcal{H}\bar{\otimes} L^2(\mathcal{A}_{\theta}))_q$.
Thus (\ref{conj:eq}) tells us that for any
$f\in (\mathcal{A}\otimes \mathcal{A}_{\theta})^{\sigma\otimes
  \tau^{-1}}$ the two restrictions $f|_{(\mathcal{H}\bar{\otimes}
    L^2(\mathcal{A}_{\theta}))_0}$ and $f|_{(\mathcal{H}\bar{\otimes}
    L^2(\mathcal{A}_{\theta}))_q}$ are unitarily conjugate to each other.
Hence
\begin{eqnarray*}
\pa f|_{(\mathcal{H}\bar{\otimes}
    L^2(\mathcal{A}_{\theta}))_0}\pa=
\pa f|_{(\mathcal{H}\bar{\otimes}
    L^2(\mathcal{A}_{\theta}))_q}\pa
\end{eqnarray*}
for all $q\in \Ze^n$. Then (\ref{same norm:eq}) follows immediately.
\end{proof}

Now Proposition~\ref{same norm theta:prop} is just a consequence of
Proposition~\ref{same norm:prop} applied to $\mathcal{A}=
C(M,Cl^{\Ce}M)$.

\section{Proof of Theorem~\ref{theta-deform metric:thm}}
\label{proof:sec}

In this section we prove Theorem~\ref{theta-deform metric:thm}
by verifying the conditions in Theorem~\ref{use of action:thm} for
the quadruple $(C(M_{\theta}), L_{\theta}, \mathbb{T}^n, I\otimes
\tau)$. Clearly $L_{\theta}$ satisfies the reality condition
(\ref{real:eq}). The condition (1) is already verified in (\ref{Lltheta to
  Ltheta:eq}).

Let $\alpha=I\otimes \tau$, and let $\hat{\alpha}=I\hat{\otimes}
\tau$. Notice that $\alpha$ is in fact an action of $\mathbb{T}^n$
on $C(M)\otimes \mathcal{A}_{\theta}$, under which $C(M_{\theta})$
is stable. For any $f\in C(M_{\theta})$ and any continuous
function $\varphi:\mathbb{T}^n\rightarrow \Ce$ clearly
$\alpha_{\varphi}(f)$ doesn't depend on whether we think of $f$ as
being in $C(M_{\theta})$ or $C(M)\otimes \mathcal{A}_{\theta}$,
where $\alpha_{\varphi}$ is the linear map on $C(M)\otimes
\mathcal{A}_{\theta}$ or $C(M_{\theta})$ defined in
Lemma~\ref{proj TVS:lemma}.

Now we verify the condition (2):

\begin{proposition} \label{nonincrease:prop}
Let $\varphi\in C(\mathbb{T}^n)$.
Then
$\Psi(C^{\infty}(M)\hat{\otimes}\mathcal{A}^{\infty}_{\theta})$
and $\Psi_{\theta}(C^{\infty}(M_{\theta}))$ are both stable under
$\alpha_{\varphi}$. We have
\begin{eqnarray} \label{nonincrease:eq}
L^D\circ \alpha_{\varphi}
\le \pa \varphi\pa_1\cdot L^D
\end{eqnarray}
on $C(M)\otimes \mathcal{A}_{\theta}$,
and
\begin{eqnarray} \label{nonincrease theta:eq}
L_{\theta}\circ \alpha_{\varphi}
\le \pa\varphi\pa_1\cdot L_{\theta}
\end{eqnarray}
on $C(M_{\theta})$.
\end{proposition}
\begin{proof}
For any $f\in
C^{\infty}(M)\hat{\otimes}\mathcal{A}^{\infty}_{\theta}$ by
Lemma~\ref{exchange integ:lemma} we have
$\alpha_{\varphi}(\Psi(f))= \Psi(\hat{\alpha}_{\varphi}(f))\\ \in
\Psi(C^{\infty}(M, S)\hat{\otimes}\mathcal{A}^{\infty}_{\theta})$.
So $\Psi(C^{\infty}(M)\hat{\otimes}\mathcal{A}^{\infty}_{\theta})$
is stable under $\alpha_{\varphi}$. For any $g\in
C^{\infty}(M_{\theta})$ by Lemma~\ref{exchange
  integ:lemma}
we have $\hat{\alpha}_{\varphi}(g)\in C^{\infty}(M_{\theta})$.
Then $\alpha_{\varphi}(\Psi_{\theta}(g))=\alpha_{\varphi}(\Psi(g))=
\Psi(\hat{\alpha}_{\varphi}(g))\in \Psi(C^{\infty}(M_{\theta}))=
\Psi_{\theta}(C^{\infty}(M_{\theta}))$.
So $\Psi_{\theta}(C^{\infty}(M_{\theta}))$
is also stable under
$\alpha_{\varphi}$.

Notice that $D^{L^2}$ is invariant under the conjugation of
$\sigma \bar{\otimes} I$, and hence $D^{L^2}_{\theta}$ is invariant under
the conjugation of the restriction of $\sigma\bar{\otimes}I$ to
$\mathcal{H}_{\theta}$. Then clearly $L^D$ and $L_{\theta}$ are
invariant under $\alpha$. Also notice that seminorms
defined by commutators are lower semicontinuous
\cite[Proposition 3.7]{Rieffel99b}.
Then (\ref{nonincrease:eq}) and
(\ref{nonincrease theta:eq}) follow
from Remark~\ref{use of action:remark}(3).
\end{proof}

We proceed to verify the conditions (3) and (4).
For
each $q\in \Ze^n=\widehat{\mathbb{T}^n}$ let $(C(M_{\theta}))_q$ be
the $q$-isotypic component of
$C(M_{\theta})$ under $\alpha$ throughout the rest of this
section. Also let $(C(M))_q$ and $(C^{\infty}(M))_q$ be the
$q$-isotypic components of $C(M)$ and $C^{\infty}(M)$ under $\sigma$.
We need:

\begin{lemma} \label{q component:lemma}
For each $q\in \Ze^n$
we have
\begin{eqnarray} \label{C q:eq}
(C(M_{\theta}))_q=(C(M))_q\otimes u_q,
\end{eqnarray}
and
\begin{eqnarray} \label{C^infty q:eq}
(C(M_{\theta}))_q\cap \Psi_{\theta}(C^{\infty}(M_{\theta}))=
(C^{\infty}(M))_q\otimes u_q.
\end{eqnarray}
\end{lemma}
\begin{proof}
Let $V=C^{\infty}(M)\hat{\otimes}\mathcal{A}^{\infty}_{\theta}$,
and let $W=C(M)\otimes\mathcal{A}_{\theta}$. Let $V_q$ and $W_q$
be the $q$-isotypic component of $V$ and $W$ under $\hat{\alpha}$
and $\alpha$ respectively. By similar arguments as in
Lemma~\ref{injective Phi:lemma}, we have $V_q=C^{\infty}(M)\otimes
u_q$ and $W_q=C(M)\otimes u_q$. Then
\begin{eqnarray*}
(C(M_{\theta}))_q=
W_q\cap W^{\sigma\otimes\tau^{-1}}
=(C(M)\otimes u_q)^{\sigma\otimes
  \tau^{-1}}=(C(M))_q\otimes u_q.
\end{eqnarray*}
Since $\Psi$ is injective, we also have
\begin{eqnarray*}
(C(M_{\theta}))_q\cap \Psi_{\theta}(C^{\infty}(M_{\theta}))
&=&
\Psi(V_q\cap V^{\sigma\hat{\otimes}\tau^{-1}})
=\Psi((C^{\infty}(M)\otimes u_q)^{\sigma\hat{\otimes}\tau^{-1}})\\
&=& \Psi((C^{\infty}(M))_q\otimes
u_q))=(C^{\infty}(M))_q\otimes u_q.
\end{eqnarray*}
\end{proof}

The geodesic distance on $M$ defines a seminorm $L_{\rho}$ on $C(M)$
via (\ref{dist to Lip:eq}).
This makes $C(M)$ into a compact quantum metric space (see the
discussion after Lemma 4.6 in \cite{Rieffel00}). Let $r_M$ be the
radius. Define a new seminorm $L$ on $C(M)$ by $L=L_{\rho}$ on
$C^{\infty}(M)$, and $L=+\infty$ on $C(M)\setminus C^{\infty}(M)$.
Since $L\ge L_{\rho}$, by Proposition~\ref{criterion of Lip:prop}
clearly $L$ is also a Lip-norm and has radius no bigger than
$r_M$. It is well known \cite{Connes89, Connes94} that
\begin{eqnarray} \label{L to d:eq}
L(f)=\pa df\pa =\pa [D, f]\pa
\end{eqnarray}
for all $f\in C^{\infty}(M)$, where we denote the closure of $D$ on
$\mathcal{H}$ also by $D$.
Notice that for any
$f=f_q\otimes u_q\in (C^{\infty}(M))_q\otimes u_q$ we have
\begin{eqnarray*}
L^D(f)=\pa [D^{L^2}, f]\pa =\pa [D, f_q]\otimes u_q\pa =
\pa [D, f_q]\pa\overset{(\ref{L to d:eq})}=L(f_q).
\end{eqnarray*}
Combining this with (\ref{Ltheta to LD:eq}), we get
\begin{eqnarray} \label{Ltheta to L:eq}
L_{\theta}(f_q\otimes u_q)=L(f_q)
\end{eqnarray}
for $f_q\otimes u_q\in (C^{\infty}(M))_q\otimes u_q$.
From (\ref{Ltheta to L:eq}), (\ref{C q:eq}) and (\ref{C^infty q:eq})
we see that $L_{\theta}$ restricted to $(C(M_{\theta}))_q$ can be identified
with $L$ restricted to $(C(M))_q$. Then conditions (3) and (4)
of Theorem~\ref{use of action:thm}
follow immediately. Then Theorem~\ref{theta-deform metric:thm} is just
a consequence
of Theorem~\ref{use of action:thm} applied to
$(C(M_{\theta}), L_{\theta}, \mathbb{T}^n, \alpha)$.
We also see that $(C(M_{\theta}), L_{\theta})$ has radius no bigger
than $r_M+C\int_{\mathbb{T}^n}\mathnormal{l}(x)\, dx$.

\end{document}